\documentclass[12pt]{article}
\usepackage{a4wide}
\usepackage[centertags,leqno]{amsmath}
\usepackage{amssymb,amsfonts}
\usepackage{color}
\usepackage{graphicx}
\usepackage[colorlinks=true]{hyperref}

\newtheorem{theorem}{Theorem}[section]

\newtheorem{corollary}[theorem]{Corollary}

\newtheorem{lemma}[theorem]{Lemma}
\newtheorem{proposition}[theorem]{Proposition}
\newtheorem{remark}[theorem]{Remark}

\definecolor{violet}{rgb}{0.5,0,0.5}
\definecolor{orange}{cmyk}{0,0.3,0.7,0}

\newtheorem{definition}{Definition}

\numberwithin{equation}{section}

\newcommand{\qed}{\rule{2mm}{2mm}}

\newcommand{\eqdef}{\stackrel{{\mathrm {def}}}{=}}

\renewcommand{\colon}{:\,}
\newcommand{\proof}{{\em Proof. }}

\newcommand{\RR}{\mathbb{R}}

\newcommand{\NN}{\mathbb{N}}

\begin{document}

\begingroup\Large\bf
\begin{center}
  New patterns of travelling waves in\\
  the generalized Fisher\--Kolmogorov equation
\end{center}
\endgroup

\vspace{0.5cm}
\begingroup\rm
\begin{center}
\centerline{Pavel {\sc Dr\'{a}bek}
\footnote{
The work of Pavel Dr\'abek was supported in part by
the Grant Agency of the Czech Republic (GA\v{C}R)
under Grant {\#}$13-00863$S.}
}
\begin{tabular}{c}
  Department of Mathematics and\\
  N.T.I.S. (Center of New Technologies for Information Society)\\
  University of West Bohemia\\
  P.O.~Box 314\\
  CZ-306 14 Plze\v{n}, Czech Republic\\
  {\it e}-mail: {\tt pdrabek@kma.zcu.cz}
\end{tabular}
\end{center}
\endgroup

\vspace{0.0cm}
\begin{center}
and
\end{center}

\vspace{0.0cm}
\begingroup\rm
\begin{center}
\centerline{Peter {\sc Tak\'{a}\v{c}}
\footnote{
The work of Peter Tak\'a\v{c} was supported in part by
Deutsche Forschungs\-gemeinschaft (DFG, Germany)
under Grant {\#} TA~213/16--1.}
}
\begin{tabular}{c}
  Institut f\"ur Mathematik\\
  Universit\"at Rostock\\
  Ulmenstra{\ss}e~69, Haus~3\\
  D-18055 Rostock, Germany\\
  {\it e}-mail: {\tt peter.takac@uni-rostock.de}
\end{tabular}
\end{center}
\endgroup

\vspace{0.3cm}

\begin{center}
\today
%December 29, 2002
\end{center}

\vspace{0.3cm}

%\newpage

\baselineskip=12pt
%%%%%%%%%%%%%%%%%%%%%%%%%%%%%%%%%%%%%%%%%%%%%%%%%%%%%%%%%%%%%%%%%%%%%%%%%%%
\noindent
\begingroup\footnotesize
{\bf {\sc Abstract.}}
We prove the existence and uniqueness of a family of travelling waves
in a degenerate (or singular) quasilinear parabolic problem
that may be regarded as a generalization of the semilinear
{\sc Fisher}-{\sc Kolmogorov}-{\sc Petrovski}-{\sc Piscounov}
equation for the advance of advantageous genes in biology.
Depending on the relation between the {\em nonlinear diffusion\/} and
the {\em nonsmooth reaction function\/},
which we quantify precisely, we investigate
the shape and asymptotic properties of travelling waves.
Our method is based on comparison results for semilinear ODEs.
\endgroup

\vfill
\par\vspace*{0.2cm}
\noindent
\begin{tabular}{ll}
{\bf Running head:}
& Travelling waves in the FKPP equation\\
\end{tabular}

\par\vspace*{0.2cm}
\noindent
\begin{tabular}{ll}
{\bf Keywords:}
& Fisher\--Kolmogorov equation, travelling waves,\\
& nonlinear diffusion, nonsmooth reaction function,\\
& comparison principle\\
\end{tabular}

\par\vspace*{0.2cm}
\noindent
\begin{tabular}{ll}
{\bf 2010 Mathematics Subject Classification:}
& Primary   35Q92, 35K92;\\
& Secondary 35K55, 35K65 \\
\end{tabular}
 
\newpage

\baselineskip=14pt
%%%%%%%%%%%%%%%%%%%%%%%%%%%%%%%%%%%%%%%%%%%%%%%%%%%%%%%%%%%%%%%%%%%%%%%%%%%
\section{Introduction}
\label{s:Intro}

The purpose of this article is to investigate
a very basic pattern formation in a reaction\--diffusion model, namely,
{\it travelling waves\/}.
The model is the favorite {\it Fisher\--KPP equation\/}
(or {\it Fisher\--Kolmogorov equation\/})
derived by {\sc R.~A.\ Fisher} \cite{Fisher} in $1937$
and first mathematically analyzed by
{\sc A.\ Kolmogorov}, {\sc I.\ Petrovski}, and {\sc N.\ Piscounov}
\cite{KPP} in the same year.
However, these original works (\cite{Fisher, KPP})
consider only {\it linear\/} diffusion and
(sufficiently) {\it smooth\/} (nonlinear) reaction.
In our present work, we allow for both,
a {\it\bfseries nonlinear\/} diffusion operator
(with a $(p-1)$-homogeneous quasilinear elliptic part, $1 < p < \infty$)
and
a {\it\bfseries nonsmooth\/} reaction function of H\"older class
$C^{0,\alpha - 1}(\RR)$ with $1 < \alpha < 2$,
\begin{equation}
\label{e:FKPP}
\left\{
\begin{aligned}
    \frac{\partial u}{\partial t}
& = \frac{\partial  }{\partial x}
    \left( d(u)
           \genfrac{|}{|}{}0{\partial u}{\partial x}^{p-2}
           \genfrac{}{}{}0{\partial u}{\partial x}
    \right) - f(u) \,,\quad (x,t)\in \RR\times \RR_+ \,,
\\
  u(x,t)&= U(x-ct) \quad\mbox{ for some speed }\, c\in \RR \,.
\end{aligned}
\right.
\end{equation}
Here, $\RR_+\eqdef [0,\infty)$, $1 < p < \infty$,
$d\colon \RR\to \RR$ is a positive continuous function, and
$f\colon \RR\to \RR$ is a continuous, but not necessarily smooth function
of a ``generalized'' Fisher\--KPP type (specified below).
For $p=2$, this equation describes
a deterministic version of a stochastic model
for the spatial spread of a favored gene in a population, suggested by
{\sc R.~A.\ Fisher} \cite{Fisher}.

Closely related situations and solution ideas have already been explored,
e.g., in
{\sc P. Dr\'abek}, {\sc R.~F. Man\'asevich}, and {\sc P. Tak\'a\v{c}}
\cite{DrabManTak} and in
{\sc Y.~Sh. Il'yasov} and {\sc P. Tak\'a\v{c}}
\cite{Ilyasov-Tak},
where existence, uniqueness, and stability of
{\it phase transition solutions\/}
in a Cahn\--Hilliard\--type model are investigated.
More precisely, we study the {\it\bfseries interaction\/}
between the (nonlinear) diffusion and the (nonsmooth) reaction;
in paticular, we explore their influence on the formation and
the shape of a {\it travelling wave\/} connecting
two stable (spatially constant) steady states.
Our main result, Theorem~\ref{thm-Main},
contains the {\it existence\/} and {\it uniqueness\/} of
such a family of travelling waves
(parametrized by a spatial shift).
This result is crucial for establishing long\--time
{\it front propagation\/} (convergence)
in non\-linear parabolic equations towards a travelling wave; see, e.g.,
{\sc D.~G.\ Aronson} and {\sc H.~F.\ Weinberger} \cite{AronWein},
{\sc P.~C.\ Fife} and {\sc J.~B.\ Mc{L}eod} \cite{Fife-McLeod}, and
{\sc F.\ Hamel} and {\sc N.\ Nadirashvili} \cite{Hamel-Nadira}
for $p=2$ (the original semilinear Fisher\--KPP equation),
{\sc E.\ Feireisl} et {\sc al.} \cite{FeHiPeTa}
for any $p\in (1,\infty)$, and
{\sc Q.\ Yi} and {\sc J.-N.\ Zhao} \cite{Yi-Zhao}
for $p\in (2,\infty)$ only.

In order to prove Theorem~\ref{thm-Main} we use
a phase plane transformation
(cf.\ {\sc J.~D.\ Murray} \cite{Murray-I}, {\S}13.2, pp.\ 440--441)
to investigate a nonlinear, first order ordinary differential equation
with an unknown parameter $c\in \RR$, see \eqref{BVP:FKPP:y(r)}.
Since this differential equation is supplemented by
homogeneous Dirichlet boundary conditions at both end\--points,
this boundary value problem is over\-determined.
We find a unique value of $c$ for which this problem
has a unique positive solution.
The nonlinearity in this differential equation does not satisfy
a local Lipschitz condition,
so, due to the lack of uniqueness of a solution,
the classical shooting method cannot be applied directly.
The novelty of our approach is to overcome this difficulty:
We take advantage of monotonicity properties of $f$ and treat
the differential equation in problem \eqref{BVP:FKPP:y(r)}
as an initial (at $-1$) or terminal (at $+1$) value problem.
We thus derive various comparison principles
that compensate for the lack of uniqueness in the shooting method.
Finally, we obtain new shapes for travelling waves by
asymptotic analysis near the end\--points $\mp 1$.

Since the seminal paper by
{\sc A.~M.\ Turing} \cite{Turing}
has appeared in $1952$, the (linear) diffusion has been known
to have a {\it destabilizing effect\/} on stable steady states
in a reaction\--diffusion equation with a smooth reaction function.
This effect leads to the formation of new, more complicated patterns,
for instance, in morphology used in {\sc Turing}'s work \cite{Turing}.
Among other things, we will demonstrate that the analogues of
{\sc Turing}'s findings do {\it not\/} apply universally
to our nonlinear problem setting.
(Our situation is similar to, but not identical with that treated in
 \cite{Turing}.)
We will determine a simple, exact relation between the constants
$p$ and $\alpha$ ($1 < p < \infty$ and $1 < \alpha < 2$),
when a destabilizing effect occurs
($1 < \alpha < p < 2$)
and when it does {\it not\/} occur
($1 < p\leq \alpha < \infty$).
Loosely speaking, this effect depends on the product
\begin{equation*}
  \mbox{ (``diffusion'') }\times \mbox{ (``smoothness'') }
  \equiv \frac{1}{p}\cdot \alpha < (\geq)\; 1 \,.
\end{equation*}
Although our results bear resemblance to
{\sc Turing}'s observations \cite{Turing}
made in a branching (bifurcation) setting,
his destabilizing effect occurs precisely for $p = \alpha = 2$.
In this case, we do {\em not\/} have any branching phenomenon
in our model in the sense that the values of our travelling wave stay
in the open interval $(-1,1)$.
In contrast, we speak of a branching phenomenon exactly when
the travelling wave attains one of the values $\mp 1$
at a (finite) spatial point $x_{\mp 1}\in \RR$.

In his monograph \cite{Murray}, {\S}14.9, pp.\ 424--430,
{\sc J.~D.\ Murray} studies
\begingroup\it
``nonexistence of stable spatial patterns for scalar equations in
one dimension with zero flux boundary conditions''.
\endgroup
In particular, the following interesting conclusion is derived
on p.~426:
\begingroup\it
Large diffusion prevents spatial patterning in
reaction diffusion mechanisms with zero flux boundary conditions.
\endgroup
Drawing an analogue of this conclusion to our setting,
we will show (Theorem~\ref{thm-Main})
that this is the case if $p\leq \alpha$.
To be more precise, we investigate {\it monotone\/} travelling waves in
the degenerate (or singular) second\--order para\-bolic problem
of a ``generalized'' Fisher\--KPP type \eqref{e:FKPP}.
For $p=2$ and $f$ smooth, extensive studies of travelling waves
can be found in
{\sc P.~C.\ Fife} and {\sc J.~B.\ Mc{L}eod} \cite{Fife-McLeod} and
{\sc J.~D.\ Murray} \cite{Murray-I}, {\S}13.2, pp.\ 439--444.

Recall that $d\colon \RR\to \RR$ is continuous and positive.
The function $f\colon \RR\to \RR$
is assumed to be continuous, such that
$f(\pm 1) = f(s_0) = 0$ for some $-1 < s_0 < 1$, together with
$f(s) > 0$ for every $s\in (-1,s_0)$,
$f(s) < 0$ for every $s\in (s_0,1)$, and
\begin{equation}
\label{int:f(u)}
  G(r)\eqdef \int_{-1}^r d(s)^{1/(p-1)}\, f(s) \,\mathrm{d}s > 0
    \quad\mbox{ whenever }\, -1 < r < 1 \,.
\end{equation}

Typical examples for the diffusion coefficient $d(s) > 0$ are
\begin{enumerate}
\renewcommand{\labelenumi}{({\bf d\arabic{enumi}})}
\item[({\bf d1})]
\makeatletter
\def\@currentlabel{{\bf d1}}\label{exam_1:d(s)}
\makeatother
$\;$
$d(s)\equiv 1$ for all $s\in \RR$,
which yields the $p$-Laplacian in eq.~\eqref{e:FKPP}, $1 < p < \infty$.
\item[({\bf d2})]
\makeatletter
\def\@currentlabel{{\bf d2}}\label{exam_2:d(s)}
\makeatother
$\;$
$p=2$ and $d(s) = \varphi'(s) > 0$ for all $s\in \RR$,
which yields the porous medium operator
$u\mapsto \frac{\partial^2}{\partial x^2} \varphi(u)$
on the right\--hand side of eq.~\eqref{e:FKPP}, where the function
$\varphi\colon \RR\to \RR$ is assumed to be
monotone increasing and continuous.
\end{enumerate}

An important special case of the reaction function $f$ is
$f(s) = F'(s)$ where
$F\colon \RR\to \RR$ is the ``generalized'' double\--well potential
$F(s)\equiv F_{\alpha}(s) = \frac{1}{2\alpha} |s^2 - 1|^{\alpha}$
for $s\in \RR$, where $\alpha\in (1,\infty)$; hence, $s_0 = 0$ and
\begin{equation}
\label{e:f_alpha}
    f(s)\equiv f_{\alpha}(s)
  = |s^2 - 1|^{\alpha - 2} (s^2 - 1) s
    \quad\mbox{ for }\, s\in \RR \,.
\end{equation}
Notice that, if $1 < \alpha < 2$ then the function $f(s)$ is only
$(\alpha - 1)$-H\"older continuous at the points $s = \pm 1$,
but certainly not Lipschitz continuous.
If the diffusion coefficient $d\colon \RR\to \RR$
(already assumed to be a positive continuous function)
is also even about zero, that is,
$d(-s) = d(s)$ for all $s\in \RR$, then \eqref{e:f_alpha} forces
\begin{equation}
\label{int:f(u):u=1}
  G(1)\eqdef \int_{-1}^1 d(s)^{1/(p-1)}\, f(s) \,\mathrm{d}s = 0
\end{equation}
in formula \eqref{int:f(u)}.

This article is organized as follows.
Our main results are collected in Section~\ref{s:Main}.
In the next section (Section~\ref{s:Prelim}),
the travelling wave problem for the Fisher\--KPP equation \eqref{e:FKPP},
the quasilinear ODE~\eqref{eq:FKPP},
is transformed into an equivalent problem for
an expression $y\colon (-1,1)\to \RR$,
the semilinear ODE~\eqref{eq:FKPP:y(r)}
with the boundary conditions \eqref{bc:FKPP:y(r)},
cf.\ problem~\eqref{BVP:FKPP:y(r)}.
The substituted unknown expression $y$ is a simple function of
the travelling wave $U(x)$ and its derivative $U'(x)$,
thus yielding a simple differential equation for
the travelling wave $U$, on one hand.
On the other hand, in eq.~\eqref{eq:FKPP:y(r)},
the unknown function $y$ depends solely on $U\in (-1,1)$
as an independent variable, i.e., $y = y(U)$.
This problem, with a monotone, non\--Lipschitzian nonlinearity,
is solved gradually in Sections
\ref{s:Exist} (existence) and
\ref{s:Unique} (uniqueness), respectively.
Finally, an important special case
(with a Fisher\--KPP\--type reaction function),
which is a slight generalization of~\eqref{e:f_alpha},
is treated in Section~\ref{s:Asympt}.

%%%%%%%%%%%%%%%%%%%%%%%%%%%%%%%%%%%%%%%%%%%%%%%%%%%%%%%%%%%%%%%%%%%%%%%%%%%
\section{Preliminaries}
\label{s:Prelim}

Assuming that the travelling wave takes the form
$u(x,t) = U(x-ct)$, $(x,t)\in \RR\times \RR_+$, with
$U\colon \RR\to \RR$ being strictly monotone decreasing and
continuously differentiable with $U'< 0$ on $\RR$, below,
we are able to find a {\it first integral\/} for
the second\--order equation for $U$:
\begin{equation}
\label{eq:FKPP}
    \frac{\mathrm{d} }{\mathrm{d}x}
    \left( d(U)
           \genfrac{|}{|}{}0{\mathrm{d}U}{\mathrm{d}x}^{p-2}
           \genfrac{}{}{}0{\mathrm{d}U}{\mathrm{d}x}
    \right)
  + c\, \frac{\mathrm{d}U}{\mathrm{d}x} - f(U)
  = 0 \,,
    \quad x\in \RR \,.
\end{equation}
Following the standard idea of phase plane transformation $(U,V)$
for the $p$-Laplacian
(cf.\ \cite[Sect.~1]{EGSanchez}),
we make the substitution
\begin{equation*}
  V\eqdef {}-
           d(U)
           \genfrac{|}{|}{}0{\mathrm{d}U}{\mathrm{d}x}^{p-2}
           \genfrac{}{}{}0{\mathrm{d}U}{\mathrm{d}x}
  > 0 \,,
\end{equation*}
whence
\begin{equation}
\label{e:dU/dx}
    \genfrac{}{}{}0{\mathrm{d}U}{\mathrm{d}x}
  = {}-
    \genfrac{(}{)}{}0{V}{d(U)}^{1/(p-1)} < 0 \,,
\end{equation}
and consequently look for $V = V(U)$ as a function of
$U\in (-1,1)$ that satisfies the following differential equation
obtained from eq.~\eqref{eq:FKPP}:
\begin{equation*}
%\label{eq:FKPP:V(U)}
{}- \frac{\mathrm{d}V}{\mathrm{d}U}\cdot \frac{\mathrm{d}U}{\mathrm{d}x}
  + c\, \frac{\mathrm{d}U}{\mathrm{d}x} - f(U)
  = 0 \,,
    \quad x\in \RR \,,
\end{equation*}
that is,
\begin{equation}
\label{eq:FKPP:V(U)}
    \frac{\mathrm{d}V}{\mathrm{d}U}\cdot
    \genfrac{(}{)}{}0{V}{d(U)}^{1/(p-1)}
  - c\, \genfrac{(}{)}{}0{V}{d(U)}^{1/(p-1)}
  - f(U) = 0 \,,
    \quad U\in (-1,1) \,.
\end{equation}
Finally, we multiply the last equation by
$d(U)^{1/(p-1)}$, make the substitution
$y\eqdef V^{p'} > 0$, where $p'= p/(p-1)\in (1,\infty)$,
and write $r$ in place of $U$, thus arriving at
\begin{equation*}
%\label{eq:FKPP:y(r)}
    \frac{1}{p'}\cdot \frac{\mathrm{d}y}{\mathrm{d}r}
  - c\, y^{1/p} - d(r)^{1/(p-1)}\, f(r)
  = 0 \,,
    \quad r\in (-1,1) \,.
\end{equation*}
This means that the unknown function
$y\colon (-1,1)\to (0,\infty)$ of $r$,
\begin{equation}
\label{e:y=V^p'}
  y = V^{p'} = d(U)^{p'}\,
      \genfrac{|}{|}{}0{\mathrm{d}U}{\mathrm{d}x}^{p}
  > 0 \,,
\end{equation}
must satisfy the following differential equation:
\begin{equation}
\label{eq:FKPP:y(r)}
    \frac{\mathrm{d}y}{\mathrm{d}r}
  = p'\left( c\, (y^{+})^{1/p} + g(r) \right) \,,
    \quad r\in (-1,1) \,,
\end{equation}
where $y^{+} = \max\{ y,\, 0\}$ and
$g(r)\eqdef d(r)^{1/(p-1)}\, f(r)$ satisfies the same hypotheses as $f$:
\begin{itemize}
\item[{}]
$g\colon \RR\to \RR$ is a continuous, but not necessarily smooth function,
such that
$g(\pm 1) = g(s_0) = 0$ for some $-1 < s_0 < 1$, together with
$g(s) > 0$ for every $s\in (-1,s_0)$,
$g(s) < 0$ for every $s\in (s_0,1)$, and \eqref{int:f(u)}, i.e.,
\begin{equation}
\label{int:g(r)}
  G(r)\eqdef \int_{-1}^r g(s) \,\mathrm{d}s > 0
    \quad\mbox{ whenever }\, -1 < r < 1 \,.
\end{equation}
\end{itemize}
Since we require that $U = U(x)$ be sufficiently smooth,
at least continuously differentiable, with
$U'(x)\to 0$ as $x\to \pm\infty$,
the function $y = y(r)$ must satisfy the boundary conditions
\begin{equation}
\label{bc:FKPP:y(r)}
  y(-1) = y(1) = 0 \,.
\end{equation}

The following remark on the value of $G(1)$ (${}\geq 0$) is in order.

%%%%%%%%%%%%%%%%%%%%%%%%%%%%%%%%%%%%%%%%%%%%%%%%%%%%%%%%%%%%%%%%%%%%%%%
%%%%%    y(r) >< G(r) is a solution (Remark)    %%%%%%%%%%%%%%%%%%%%%%%
%%%%%%%%%%%%%%%%%%%%%%%%%%%%%%%%%%%%%%%%%%%%%%%%%%%%%%%%%%%%%%%%%%%%%%%
\begin{remark}\label{rem-sol:G(r)}\nopagebreak
\begingroup\rm
Since the integrand
$s\mapsto g(s)\eqdef d(s)^{1/(p-1)}\, f(s)\colon (-1,1)\to \RR$
in the function $G(r)$,
defined in~\eqref{int:f(u)} for $r\in (-1,1)$,
is continuous and absolutely integrable over $(-1,1)$,
we conclude that
$G\colon [-1,1]\to \RR$ is absolutely continuous.
In particular, ineq.~\eqref{int:g(r)} forces $G(1)\geq 0$.
We will see later that the case $G(1) = 0$ guarantees
the existence of a {\em stationary solution\/}
to problem~\eqref{e:FKPP}, i.e., $c=0$, whereas
the case $G(1) > 0$ renders a {\em travelling wave\/}, i.e., $c\neq 0$;
more precisely, with $c<0$, cf.\
Theorem~\ref{thm-Main} below.
Indeed, both,
the stationary solution (for $c = 0$) and
the travelling wave (for $c < 0$)
will be obtained from eq.~\eqref{eq:FKPP:y(r)}
by means of the transformation defined by eqs.\
\eqref{e:dU/dx} and \eqref{e:y=V^p'}.
\endgroup
\end{remark}
%%%%%%%%%%%%%%%%%%%%%%%%%%%%%%%%%%%%%%%%%%%%%%%%%%%%%%%%%%%%%%%%%%%%%%%
\par\vskip 10pt

%%%%%%%%%%%%%%%%%%%%%%%%%%%%%%%%%%%%%%%%%%%%%%%%%%%%%%%%%%%%%%%%%%%%%%%%%%%
\section{Main Results}
\label{s:Main}

We assume that
$d\colon \RR\to \RR$ is a positive continuous function.
In fact, we need only
$d\colon [-1,1]\to \RR$ to be continuous and positive.
Let us recall that the function
$g(r) = d(r)^{1/(p-1)}\, f(r)$ of $r\in [-1,1]$
is assumed to satisfy the hypotheses formulated in the previous section
(Section~\ref{s:Prelim})
before Remark~\ref{rem-sol:G(r)},
in particular, ineq.~\eqref{int:g(r)}.

Our main result concerning the {\em travelling waves\/}
$u(x,t) = U(x-ct)$, $(x,t)\in \RR\times \RR_+$, $c\in \RR$, for
problem~\eqref{e:FKPP} is as follows.

%%%%%%%%%%%%%%%%%%%%%%%%%%%%%%%%%%%%%%%%%%%%%%%%%%%%%%%%%%%%%%%%%%%%%%%
%%%%%    Main Results ($\alpha / p$) (Theorem)    %%%%%%%%%%%%%%%%%%%%%
%%%%%%%%%%%%%%%%%%%%%%%%%%%%%%%%%%%%%%%%%%%%%%%%%%%%%%%%%%%%%%%%%%%%%%%
\begin{theorem}\label{thm-Main}
Let\/ $G(1) > 0$.
Then there exists a unique number\/ $c^{\ast}\in \RR$ such that\/
problem~\eqref{e:FKPP} with $c = c^{\ast}$
possesses a travelling wave solution
$u(x,t) = U(x - c^{\ast} t)$, $(x,t)\in \RR\times \RR_+$,
where $U = U(\xi)$ is
a monotone decreasing and continuously differentiable function on $\RR$
taking values in $[-1,1]$.
Furthermore, we have $c^{\ast} < 0$ and the set\/
\begin{math}
  \{ \xi\in \RR\colon U'(\xi) < 0\}
\end{math}
is a nonempty open interval $(x_1, x_{-1})\subset \RR$ with
\begin{equation}
\label{e:U(x),x_+-1}
  \lim_{\xi\to (x_1)+}    U(\xi) =  1 \qquad\mbox{ and }\qquad
  \lim_{\xi\to (x_{-1})-} U(\xi) = -1 \,.
\end{equation}
\end{theorem}
%%%%%%%%%%%%%%%%%%%%%%%%%%%%%%%%%%%%%%%%%%%%%%%%%%%%%%%%%%%%%%%%%%%%%%%
%\par\vskip 10pt

\par\vskip 10pt
%%%%%%%%%%%%%%%%%%%%%%%%%%%%
\proof
The existence and uniqueness of $c^{\ast}\in \RR$ follow from
the transformation of eq.~\eqref{eq:FKPP}
into eq.~\eqref{eq:FKPP:y(r)}
with the boundary conditions \eqref{bc:FKPP:y(r)}
and subsequent application of Theorem~\ref{thm-FKPP:y(r)} below
to the boundary value problem
\eqref{eq:FKPP:y(r)}, \eqref{bc:FKPP:y(r)}.
We get also $c^{\ast} < 0$.

Inserting $y = y_{ c^{\ast} }\colon (-1,1)\to \RR$
into eq.~\eqref{e:y=V^p'} we obtain
\begin{equation*}
\left\{
\begin{aligned}
& V(r) = y_{ c^{\ast} }(r)^{1/p'} = y_{ c^{\ast} }(r)^{(p-1)/p} > 0
  \quad\mbox{ for }\, r\in (-1,1) \,;
\\
& V(-1) = V(1) = 0 \,.
\end{aligned}
\right.
\end{equation*}
We insert $V = V(U)$ as a function of
$U\in (-1,1)$ into eq.~\eqref{e:dU/dx}, thus arriving at
\begin{equation*}
%\label{e:dU/dx(U)}
    \genfrac{}{}{}0{\mathrm{d}U}{\mathrm{d}x}
  = {}-
    \genfrac{(}{)}{}0{V(U)}{d(U)}^{1/(p-1)} < 0
  \quad\mbox{ for }\, U\in (-1,1) \,.
\end{equation*}
Separation of variables above yields
\begin{equation*}
%\label{e:dx=-dU/U}
  \mathrm{d}x
  = {}-
    \genfrac{(}{)}{}0{d(U)}{V(U)}^{1/(p-1)} \,\mathrm{d}U
  \quad\mbox{ for }\, U\in (-1,1) \,.
\end{equation*}
Finally, we integrate the last equation to arrive at
\begin{equation}
\label{int:dx=-dU/U}
    x(U)
  = x(0) - \int_0^U
    \genfrac{(}{)}{}0{d(r)}{V(r)}^{1/(p-1)} \,\mathrm{d}r
  \quad\mbox{ for }\, U\in (-1,1) \,,
\end{equation}
where $x(0) = x_0\in \RR$ is an arbitrary constant.
We remark that $V(r)\to 0$ as $r\to \pm 1$.

Consequently, both monotone limits below exist,
\begin{equation}
\label{e:x(U),U=+-1}
  x_{-1}\eqdef \lim_{U\to (-1)+} x(U)
    \qquad\mbox{ and }\qquad
  x_1   \eqdef \lim_{U\to (+1)-} x(U) \,,
\end{equation}
and satisfy
$-\infty\leq x_1 < x_0 < x_{-1}\leq +\infty$.
Thus, the function
$x\colon (-1,1)\to (x_1, x_{-1})\subset \RR$
is a diffeomorphism of the open interval $(-1,1)$ onto $(x_1, x_{-1})$
satisfying
$\frac{\mathrm{d}x}{\mathrm{d}U} < 0$ in $(-1,1)$
together with \eqref{e:x(U),U=+-1}.
This implies the remaining part of our theorem, especially
the limits in \eqref{e:U(x),x_+-1}.
\qed
%%%%%%%%%%%%%%%%%%%%%%%%%%%%
\par\vskip 10pt

It depends on the asymptotic behavior of the function $g = g(r)$,
$r\in (-1,1)$, near the points $\pm 1$ whether
$U'< 0$ holds on the entire real line $\RR = (-\infty, +\infty)$
or else
$U'< 0$ on a nonempty open interval $(x_1, x_{-1})\subset \RR$ and
$U'= 0$ on its nonempty complement $\RR\setminus (x_1, x_{-1})$ with
$x_1 > -\infty$ and/or $x_{-1} < +\infty$.
More precisely, we assume that
$g\colon \RR\to \RR$ has the following asymptotic behavior near $\pm 1$:

There are constants
$\gamma^{\pm}, \gamma_0^{\pm}\in (0,\infty)$ such that
\begin{equation}
\label{e:g(r),r=+-1}
    \lim_{r\to 1-}  \frac{g(r)}{ (1-r)^{\gamma^{+}} }
  = {}- \gamma_0^{+} \qquad\mbox{ and }\qquad
    \lim_{r\to -1+} \frac{g(r)}{ (1+r)^{\gamma^{-}} }
  = \gamma_0^{-} \,.
\end{equation}

For instance, if
$g(s)\eqdef d(s)^{1/(p-1)}\, f(s)$ for $s\in (-1,1)$
is as in Remark~\ref{rem-sol:G(r)} and
$f(s)\equiv f_{\alpha}(s)$ is defined by \eqref{e:f_alpha}, then we have
$\alpha = 1 + \gamma^{\pm}$.

We have the following conclusions for the limits
$x_{\mp 1}\eqdef \lim_{ U\to (\mp 1)\pm } x(U)$, where
$x = x(U)$ is a solution of
\begin{equation}
\label{e:dx/dU}
    \genfrac{}{}{}0{\mathrm{d}x}{\mathrm{d}U}
  = {}-
    \genfrac{(}{)}{}0{d(U)}{V(U)}^{1/(p-1)} < 0
    \quad\mbox{ for }\, -1 < U < 1 \,,
\end{equation}
cf.\ eq.~\eqref{e:dU/dx},
given explicitely by formula~\eqref{e:x(U),U=+-1}.
Notice that, in this notation,
$x\colon (-1,1)\to (x_1, x_{-1})$
is the inverse function of the restriction of $U$ to the interval
$(x_1, x_{-1})$ in which $U'< 0$, i.e.,
$x(r) = U^{-1}(r)$ for $-1 < r < 1$.

%%%%%%%%%%%%%%%%%%%%%%%%%%%%%%%%%%%%%%%%%%%%%%%%%%%%%%%%%%%%%%%%%%%%%%%
%%%%%    Main Result for $x_{\pm 1}$ (Theorem)    %%%%%%%%%%%%%%%%%%%%%
%%%%%%%%%%%%%%%%%%%%%%%%%%%%%%%%%%%%%%%%%%%%%%%%%%%%%%%%%%%%%%%%%%%%%%%
\begin{theorem}\label{thm-x_+-1}
Assume that the limits in eq.~\eqref{e:g(r),r=+-1} hold.  Then,

{\rm (i)}$\;$
in case $1 < p\leq 2$ we have $x_1 = -\infty$ $(x_{-1} = +\infty)$
if and only if\/
$\gamma^{+}\geq p-1$ ($\gamma^{-}\geq p-1$, respectively); whereas

{\rm (ii)}$\;$
in case $2 < p < \infty$ we have $x_1 = -\infty$ $(x_{-1} = +\infty)$
if\/
$\gamma^{+}\geq p-1$ ($\gamma^{-}\geq p-1$, respectively).
\end{theorem}
%%%%%%%%%%%%%%%%%%%%%%%%%%%%%%%%%%%%%%%%%%%%%%%%%%%%%%%%%%%%%%%%%%%%%%%
\par\vskip 10pt

This theorem is an easy combination of Theorem~\ref{thm-Main} above
with Corollary~\ref{cor-power<g(r)}
and Remark~\ref{rem-power<g(r)} in Section~\ref{s:Asympt}.

The following corollary to Theorem~\ref{thm-x_+-1}
for the linear diffusion case ($p=2$) is obvious:

%%%%%%%%%%%%%%%%%%%%%%%%%%%%%%%%%%%%%%%%%%%%%%%%%%%%%%%%%%%%%%%%%%%%%%%
%%%%%    Main Result for $x_{\pm 1}$ (Corollary)    %%%%%%%%%%%%%%%%%%%
%%%%%%%%%%%%%%%%%%%%%%%%%%%%%%%%%%%%%%%%%%%%%%%%%%%%%%%%%%%%%%%%%%%%%%%
\begin{corollary}\label{cor-x_+-1}
Let\/ $p=2$.  Then
$x_1$ $(x_{-1})$ is finite if\/
$\gamma^{+} < 1$ ($\gamma^{-} < 1$, respectively) and\/
$x_1 = -\infty$ $(x_{-1} = +\infty)$
if\/
$\gamma^{+}\geq 1$ ($\gamma^{-}\geq 1$, respectively).
\end{corollary}
%%%%%%%%%%%%%%%%%%%%%%%%%%%%%%%%%%%%%%%%%%%%%%%%%%%%%%%%%%%%%%%%%%%%%%%
\par\vskip 10pt

The results stated above hinge upon
the existence and uniqueness results for
eq.~\eqref{eq:FKPP:y(r)} with the Dirichlet boundary conditions
\eqref{bc:FKPP:y(r)}, i.e.,
for the following boundary value problem:
\begin{equation}
\label{BVP:FKPP:y(r)}
\left\{
\begin{aligned}
&   \frac{\mathrm{d}y}{\mathrm{d}r}
  = p'\left( c\, (y^{+})^{1/p} + g(r) \right) \,,
    \quad r\in (-1,1) \,;
\\
& y(-1) = y(1) = 0 \,,
\end{aligned}
\right.
\end{equation}
with the parameter $c\in \RR$ to be determined.
We recall that eq.~\eqref{eq:FKPP:y(r)}
has been obtained from eq.~\eqref{e:dU/dx}
by means of the substitution in~\eqref{e:y=V^p'}.

The following result for problem~\eqref{BVP:FKPP:y(r)}
is of independent interest
(cf.\ \cite{FeHiPeTa, Yi-Zhao}).

%%%%%%%%%%%%%%%%%%%%%%%%%%%%%%%%%%%%%%%%%%%%%%%%%%%%%%%%%%%%%%%%%%%%%%%
%%%%%    Main Result - auxiliary equation (Theorem)    %%%%%%%%%%%%%%%%
%%%%%%%%%%%%%%%%%%%%%%%%%%%%%%%%%%%%%%%%%%%%%%%%%%%%%%%%%%%%%%%%%%%%%%%
\begin{theorem}\label{thm-FKPP:y(r)}
Let\/ $G(1) > 0$.
Then there exists a unique number $c^{\ast}\in \RR$ such that\/
problem~\eqref{BVP:FKPP:y(r)} with $c = c^{\ast}$ has a unique solution
$y = y_{ c^{\ast} }\colon (-1,1)\to \RR$.
Furthermore, we have $c^{\ast} < 0$ and\/
$y_{ c^{\ast} } > 0$ in $(-1,1)$.
\end{theorem}
%%%%%%%%%%%%%%%%%%%%%%%%%%%%%%%%%%%%%%%%%%%%%%%%%%%%%%%%%%%%%%%%%%%%%%%
%\par\vskip 10pt

\par\vskip 10pt
%%%%%%%%%%%%%%%%%%%%%%%%%%%%
\proof
The existence and uniqueness of $c^{\ast}\in \RR$ follow from
Corollary~\ref{cor-uniq_c^*}
(with the existence established before in Corollary~\ref{cor-c<c^*}).
We get also $c^{\ast} < 0$.
\qed
%%%%%%%%%%%%%%%%%%%%%%%%%%%%
\par\vskip 10pt

%%%%%%%%%%%%%%%%%%%%%%%%%%%%%%%%%%%%%%%%%%%%%%%%%%%%%%%%%%%%%%%%%%%%%%%%%%%
\section{Existence Result for $c^{\ast}$}
\label{s:Exist}

In order to verify that only the case $c<0$ can yield a travelling wave
for problem~\eqref{e:FKPP},
we prove the following simple lemma for
the boundary value problem~\eqref{BVP:FKPP:y(r)}.

%%%%%%%%%%%%%%%%%%%%%%%%%%%%%%%%%%%%%%%%%%%%%%%%%%%%%%%%%%%%%%%%%%%%%%%
%%%%%    y(u) >< G(u) is a solution (Lemma)    %%%%%%%%%%%%%%%%%%%%%%%%
%%%%%%%%%%%%%%%%%%%%%%%%%%%%%%%%%%%%%%%%%%%%%%%%%%%%%%%%%%%%%%%%%%%%%%%
\begin{lemma}\label{lem-sol:G(u)}
Let\/ $c\in \RR$ and $y_0(r) = p'\, G(r)$ for\/ $-1\leq r\leq 1$.

{\rm (i)}$\;$
Let\/ $c = 0$.
Then the function $y_0(r) = p'\, G(r)$ of\/ $r\in [-1,1]$
is a solution to problem~\eqref{BVP:FKPP:y(r)}
if and only if $G(1) = 0$.

{\rm (ii)}$\;$
If\/ $c>0$ $($$c\leq 0$, respectively$)$ then every solution
$y\colon [-1,1]\to \RR$
to the initial value problem for eq.~\eqref{eq:FKPP:y(r)}
with the initial condition $y(-1) = 0$ satisfies
$y > y_0$ $($$y\leq y_0$$)$ throughout $(-1,1]$.
\end{lemma}
%%%%%%%%%%%%%%%%%%%%%%%%%%%%%%%%%%%%%%%%%%%%%%%%%%%%%%%%%%%%%%%%%%%%%%%
%\par\vskip 10pt

\par\vskip 10pt
%%%%%%%%%%%%%%%%%%%%%%%%%%%%
\proof
Part~{\rm (i)} is a trivial consequence of \eqref{int:g(r)}
combined with the properties of~$g$.

Part~{\rm (ii)}:
First, let $c>0$.
Then, clearly,
$y'\geq y_0' = p'\, g$ throughout $(-1,1)$, together with
$y'\not\equiv y_0'$ in $(-1,r)$ for every $r\in (-1,1)$, which forces
$y > y_0$ throughout $(-1,1]$, by a simple integration of
eq.~\eqref{eq:FKPP:y(r)} over the interval $[-1,r]$.

The case $c\leq 0$ is analogous, by reversing
the (nonstrict) inequalities.
\qed
%%%%%%%%%%%%%%%%%%%%%%%%%%%%
\par\vskip 10pt

For $c\geq 0$, Lemma~\ref{lem-sol:G(u)}
has the following obvious consequence.

%%%%%%%%%%%%%%%%%%%%%%%%%%%%%%%%%%%%%%%%%%%%%%%%%%%%%%%%%%%%%%%%%%%%%%%
%%%%%    y(u) >< G(u) is a solution (Corollary)    %%%%%%%%%%%%%%%%%%%%
%%%%%%%%%%%%%%%%%%%%%%%%%%%%%%%%%%%%%%%%%%%%%%%%%%%%%%%%%%%%%%%%%%%%%%%
\begin{corollary}\label{cor-sol:G(u)}
If\/ {\rm either\/} {\rm (i)}  $c=0$ and\/ $G(1) > 0$
{\rm or else\/}     {\rm (ii)} $c>0$ and\/ $G(1)\geq 0$,
then eq.~\eqref{eq:FKPP:y(r)} possesses {\rm no\/} solution
$y\colon [-1,1]\to \RR$
satisfying the boundary conditions $y(-1) = 0$ and\/ $y(1)\leq 0$.
\end{corollary}
%%%%%%%%%%%%%%%%%%%%%%%%%%%%%%%%%%%%%%%%%%%%%%%%%%%%%%%%%%%%%%%%%%%%%%%
\par\vskip 10pt

We conclude that it suffices to investigate the case $c<0$
for finding a travelling wave to problem~\eqref{e:FKPP}.

%%%%%%%%%%%%%%%%%%%%%%%%%%%%%%%%%%%%%%%%%%%%%%%%%%%%%%%%%%%%%%%%%%%%%%%
%%%%%    y(u) >< G(u) is a solution (Remark)    %%%%%%%%%%%%%%%%%%%%%%%
%%%%%%%%%%%%%%%%%%%%%%%%%%%%%%%%%%%%%%%%%%%%%%%%%%%%%%%%%%%%%%%%%%%%%%%
\begin{remark}\label{rem-sol:G(u)}\nopagebreak
\begingroup\rm
Given any $c\leq 0$, the nonlinearity
$y\mapsto c\, (y^{+})^{1/p}\colon \RR\to \RR$
is nonincreasing and, thus, satisfies a one\--sided Lipschitz condition,
which guarantees that eq.~\eqref{eq:FKPP:y(r)}
with the initial condition $y(-1) = 0$ possesses
a {\it unique\/} solution
$y\colon [-1,1]\to \RR$, by two well\--known
existence and uniqueness theorems from
{\sc Ph.~Hartman}'s monograph \cite{Hartman},
Theorem 2.1, p.~10 (existence due to Peano), and
Theorem 6.2, p.~34 (one\--sided uniqueness), respectively.
This solution depends continuously on $c\leq 0$, by
\cite[Theorem 2.1, p.~94]{Hartman},
$y\equiv y_c$, and it may change sign.
Furthermore, for every fixed $r\in (-1,1]$, the function
$c\mapsto y_c(r)\colon \RR_-\to \RR$
is monotone increasing, by
\cite[Corollary 4.2, p.~27]{Hartman}.
Here, we have denoted $\RR_-\eqdef (-\infty,0]$.
We investigate this dependence after this remark.

Even if we do not need to exploit the case $c>0$ any more,
we still need to consider
the terminal value problem for eq.~\eqref{eq:FKPP:y(r)}
with the terminal condition $y(1) = 0$, i.e.,
the ``backward'' initial value problem for eq.~\eqref{eq:FKPP:y(r)}
with $y(1) = 0$.
Unfortunately, if $c<0$, we do not have a uniqueness result of
the kind described above for this problem.
To overcome this difficulty, in the next section
(Section~\ref{s:Unique})
we establish a comparison result in
Proposition~\ref{prop-y(1):sub/sup}{\rm (b)}.
\endgroup
\end{remark}
%%%%%%%%%%%%%%%%%%%%%%%%%%%%%%%%%%%%%%%%%%%%%%%%%%%%%%%%%%%%%%%%%%%%%%%
\par\vskip 10pt

More detailed results concerning
existence, uniqueness, monotone dependence on certain parameters, and
other qualitative properties of solutions to
the initial value problem for eq.~\eqref{eq:FKPP:y(r)}
with the initial condition $y(-1) = 0$ will be established
in Corollary~\ref{cor-y_c:monot}.

In view of
Remark~\ref{rem-sol:G(r)} and Corollary~\ref{cor-sol:G(u)} above,
from now on we assume
\begin{equation}
\label{int>0:f(u):u=1}
  G(1)\eqdef \int_{-1}^1 d(s)^{1/(p-1)}\, f(s) \,\mathrm{d}s > 0
\end{equation}
in formula \eqref{int:f(u)}.
Equivalently, we have
$y_0(1) = p'\, G(1) > 0$.
We will show that there is a unique constant $c < 0$ such that
the (unique) solution $y_c\colon [-1,1]\to \RR$
to the initial value problem for eq.~\eqref{eq:FKPP:y(r)}
with the initial condition $y_c(-1) = 0$ satisfies also
the terminal condition $y_c(1) = 0$.
To find this constant, let us define
\begin{equation}
\label{def:c^*}
  c^{\ast}\eqdef \inf\left\{ c\in \RR\colon y_c(1) > 0\right\} \,.
\end{equation}
As expected, we will show that precisely $c^{\ast}$
is the desired value of the constant~$c$; see
Corollaries \ref{cor-c<c^*} and~\ref{cor-uniq_c^*}.

We need the following two technical lemmas.

%%%%%%%%%%%%%%%%%%%%%%%%%%%%%%%%%%%%%%%%%%%%%%%%%%%%%%%%%%%%%%%%%%%%%%%
%%%%%    y_c(r) > 0 for -1 < r < -1 + delta (Lemma)    %%%%%%%%%%%%%%%%
%%%%%%%%%%%%%%%%%%%%%%%%%%%%%%%%%%%%%%%%%%%%%%%%%%%%%%%%%%%%%%%%%%%%%%%
\begin{lemma}\label{lem-y_c>0}
Let\/ $c\in \RR$ be arbitrary.

{\rm (i)}$\;$
There is some
$\delta\equiv \delta(c)\in (0,2)$
such that $y_c(r) > 0$ holds for all\/ $r\in (-1, -1+\delta)$.

{\rm (ii)}$\;$
If\/ $y_c(1)\geq 0$
then we have also $y_c(r) > 0$ for all\/ $r\in (-1,1)$.
\end{lemma}
%%%%%%%%%%%%%%%%%%%%%%%%%%%%%%%%%%%%%%%%%%%%%%%%%%%%%%%%%%%%%%%%%%%%%%%
%\par\vskip 10pt

\par\vskip 10pt
%%%%%%%%%%%%%%%%%%%%%%%%%%%%
\proof
Part~{\rm (i)}:
On the contrary, if no such number $\delta\in (0,2)$ exists,
then there is a sequence
$\{ r_n\}_{n=1}^{\infty} \subset (-1,1)$ such that
$r_n\searrow -1$ as $n\nearrow \infty$, and
$y_c(r_n)\leq 0$ for all $n\in \NN$.
Given any fixed $n\in \NN$, we cannot have
$y_c(s)\leq 0$ for all $s\in (-1,r_n]$ since, otherwise,
by eq.~\eqref{eq:FKPP:y(r)}, the function $y_c$ would satisfy
\begin{equation}
\label{eq:FKPP:y_c(r)}
  \frac{\mathrm{d}y_c}{\mathrm{d}r} = p'\, g(r) \,,
    \quad r\in (-1,r_n) \,,
\end{equation}
that is, $y_c = y_0 > 0$ throughout $[-1,r_n]$, by a simple integration of
eq.~\eqref{eq:FKPP:y_c(r)} over the interval $[-1,r]$.
This conclusion would then contradict our choice $y_c(r_n)\leq 0$.
Consequently, there is a point $r_n'\in (-1,r_n)$ such that
$y_c(r_n') > 0$.
It follows that there is another point $r_n''\in (r_n',r_n]$ such that
$y_c(r_n'')\leq 0$ and $y_c'(r_n'')\leq 0$.
Since $r_n\searrow -1$ as $n\nearrow \infty$, we have also
$r_n''\to -1$ as $n\to \infty$.
But now, substituting $r = r_n''$ in eq.~\eqref{eq:FKPP:y(r)},
we arrive at $g(r_n'')\leq 0$ for every $n\in \NN$.
A contradiction with our hypothesis on $g$, namely,
$g(s) > 0$ for every $s\in (-1,s_0)$,
is obtained for every $n\in \NN$ sufficiently large, such that
$r_n''\in (-1,s_0)$.

We have proved the existence of the number $\delta$.

Part~{\rm (ii)}:
On the contrary, let us assume that the set
$\{ r\in (-1,1)\colon y_c(r)\leq 0\}$ is nonempty.
We denote by $\overline{r}\in (-1,1)$
the smallest number $r\in (-1,1)$ such that $y_c(r) = 0$.
Hence, $-1 < -1+\delta\leq \overline{r} < 1$ and
$y_c(r) > 0$ holds for every $r\in (-1,\overline{r})$.
We have
$y_c(\overline{r}) = 0$ and, therefore,
\begin{math}
  p'\, g(\overline{r}) = y_c'(\overline{r})\leq 0 \,,
\end{math}
by eq.~\eqref{eq:FKPP:y(r)}.
Our hypothesis on the nodal point $s_0$ of the function $g$ thus forces
$\overline{r}\geq s_0$.
Hence, we must have
$y_c'(r) = p'\, g(r) < 0$ for all
$r\in (\overline{r}, 1)\subset (s_0,1)$.
We conclude that $y_c(r) < 0$ holds for all
$r\in (\overline{r}, 1]$,
which contradicts our hypothesis $y_c(1)\geq 0$.

We have proved that $y_c(r) > 0$ holds for all $r\in (-1,1)$ as desired.
\qed
%%%%%%%%%%%%%%%%%%%%%%%%%%%%
\par\vskip 10pt

%%%%%%%%%%%%%%%%%%%%%%%%%%%%%%%%%%%%%%%%%%%%%%%%%%%%%%%%%%%%%%%%%%%%%%%
%%%%%    y(r) < 0 for r_0 < r < 1 and r_0\geq s_0 (Lemma)    %%%%%%%%%%
%%%%%%%%%%%%%%%%%%%%%%%%%%%%%%%%%%%%%%%%%%%%%%%%%%%%%%%%%%%%%%%%%%%%%%%
\begin{lemma}\label{lem-r_0>s_0}
Assume that\/ $c < 0$ and\/ $G(1) > 0$.
Let\/ $y\colon [-1,1]\to \RR$ be a solution to eq.~\eqref{eq:FKPP:y(r)}
satisfying\/ $y(-1) = y(r_0) = 0$ for some $r_0\in (-1,1)$ and\/
$y(r) > 0$ for all\/ $r\in (-1,r_0)$.
Then we have $r_0\geq s_0$ and\/
$y(r) < 0$ holds for all\/ $r\in (r_0,1]$.
\end{lemma}
%%%%%%%%%%%%%%%%%%%%%%%%%%%%%%%%%%%%%%%%%%%%%%%%%%%%%%%%%%%%%%%%%%%%%%%
%\par\vskip 10pt

\par\vskip 10pt
%%%%%%%%%%%%%%%%%%%%%%%%%%%%
\proof
On the contrary, suppose that $-1 < r_0 < s_0$.
Then $y'(r_0) = p'\, g(r_0) > 0$ holds by our hypotheses on $g$.
However, this is contradicts our hypothesis on $y$ requiring
$y(r) > 0 = y(r_0)$ for all $r\in (-1,r_0)$ and, thus, forcing
$y'(r_0)\leq 0$.
This proves $r_0\geq s_0$.

Consequently, we have $g(r) < 0$ for all $r\in (r_0,1)$
which entails
$y'(r) = p'\, g(r) < 0$ for all $r\in (r_0,1)$.
Integration over the interval $[r_0,r]$ yields the desired result, i.e.,
$y(r) < 0$ holds for all\/ $r\in (r_0,1]$.
\qed
%%%%%%%%%%%%%%%%%%%%%%%%%%%%
\par\vskip 10pt

%%%%%%%%%%%%%%%%%%%%%%%%%%%%%%%%%%%%%%%%%%%%%%%%%%%%%%%%%%%%%%%%%%%%%%%
%%%%%    y(u) >< G(u) is a solution (Proposition)    %%%%%%%%%%%%%%%%%%
%%%%%%%%%%%%%%%%%%%%%%%%%%%%%%%%%%%%%%%%%%%%%%%%%%%%%%%%%%%%%%%%%%%%%%%
\begin{proposition}\label{prop-c^*}
We have\/
$-\infty < c^{\ast} < 0$ and\/ $y_{c^{\ast}}(1) = 0$.
\end{proposition}
%%%%%%%%%%%%%%%%%%%%%%%%%%%%%%%%%%%%%%%%%%%%%%%%%%%%%%%%%%%%%%%%%%%%%%%
%\par\vskip 10pt

\par\vskip 10pt
%%%%%%%%%%%%%%%%%%%%%%%%%%%%
\proof
First, let us recall that the terminal value
$y_c(1)\in \RR$ is
a monotone increasing, continuous function
$c\mapsto y_c(1)\colon \RR_-\to \RR$
of the parameter $c\in \RR_-$.
These properties yield immediately
$-\infty\leq c^{\ast} < 0$ and, if $c^{\ast} > -\infty$ then also
$y_{c^{\ast}}(1) = 0$.

Thus, it remains to prove $c^{\ast} > -\infty$.
By contradiction, let us suppose that there is a sequence of numbers
$c_n\in (-\infty,0)$, such that
$y_{c_n}(1) > 0$ for every $n=1,2,3,\dots$ and
$c_n\searrow -\infty$ as $n\nearrow \infty$.
Then $y_{c_n}(r) > 0$ must hold for all $r\in (-1,1]$,
by Lemma~\ref{lem-y_c>0}, Part~{\rm (ii)}.

Next, let $z_n\colon [s_0,1]\to \RR$ be the (unique) solution of
the initial value problem
\begin{equation}
\label{eq:FKPP:z(r)}
    \frac{\mathrm{d}z_n}{\mathrm{d}r}
  = p' c_n\, (z_n^{+})^{1/p} \,,\quad r\in (s_0,1) \,;\qquad
    z_n(s_0) = y_{c_n}(s_0) \,.
\end{equation}
The monotone (increasing) dependence of the solution
$y_c\colon [-1,1]\to \RR$
on the right\--hand side of eq.~\eqref{eq:FKPP:y(r)}
(essentially due to {\sc E.\ Kamke}; see, e.g.,
 {\sc M.~W.\ Hirsch} \cite[p.~425]{Hirsch} or
 {\sc W.\ Walter}    \cite[Chapt.\ III, {\S}10]{Walter},
 Comparison Theorem on p.~112)
guarantees the following comparison result:
$0 < y_{c_n}(r)\leq z_n(r)$ for all $r\in [s_0,1]$.
However, separating variables in eq.~\eqref{eq:FKPP:z(r)}
and integrating, we obtain
\begin{equation*}
  (z_n(r))^{1/p'} = (z_n(s_0))^{1/p'} + c_n (r-s_0)
  \quad\mbox{ for all }\, r\in [s_0,1] \,.
\end{equation*}
Consequently, by the comparison result, we have also
\begin{equation*}
  0 < (y_{c_n}(r))^{1/p'}
  \leq (y_{c_n}(s_0))^{1/p'} + c_n (r-s_0)
  \quad\mbox{ for all }\, r\in [s_0,1] \,.
\end{equation*}
Setting $r=1$ and recalling $y_{c_n}(1) > 0$, we observe that
$(y_{c_n}(s_0))^{1/p'} > - c_n (1-s_0)$
holds for every $n=1,2,3,\dots$.
Since also $y_{c_n}\leq y_0$ holds throughout the interval $[-1,1]$,
by Lemma~\ref{lem-sol:G(u)}, Part~{\rm (ii)}, we conclude that
$- c_n (1-s_0) < (y_0(s_0))^{1/p'}$,
a contradiction to $-c_n\nearrow \infty$ as $n\nearrow \infty$.
We have proved $c^{\ast} > -\infty$ as desired.
\qed
%%%%%%%%%%%%%%%%%%%%%%%%%%%%
\par\vskip 10pt

Lemma~\ref{lem-r_0>s_0} and Proposition~\ref{prop-c^*}
entail the following obvious corollary.

%%%%%%%%%%%%%%%%%%%%%%%%%%%%%%%%%%%%%%%%%%%%%%%%%%%%%%%%%%%%%%%%%%%%%%%
%%%%%    y(u) >< G(u) is a solution (Corollary)    %%%%%%%%%%%%%%%%%%%%
%%%%%%%%%%%%%%%%%%%%%%%%%%%%%%%%%%%%%%%%%%%%%%%%%%%%%%%%%%%%%%%%%%%%%%%
\begin{corollary}\label{cor-c<c^*}
There is some number\/ $c\in \RR$, say, $c = c^{\ast}$ $(< 0)$,
such that problem~\eqref{BVP:FKPP:y(r)} possesses a solution
$y\colon [-1,1]\to \RR$.
This solution is unique, given by $y = y_{c^{\ast}}$, and it satisfies
$y(r) > 0$ for all\/ $r\in (-1,1)$.
\end{corollary}
%%%%%%%%%%%%%%%%%%%%%%%%%%%%%%%%%%%%%%%%%%%%%%%%%%%%%%%%%%%%%%%%%%%%%%%
\par\vskip 10pt

For this particular speed $c = c^{\ast}$,
the uniqueness of $y_{c^{\ast}}$ has been discussed
in Remark~\ref{rem-sol:G(u)}.
The sole existence of $y_{c^{\ast}}$ can be found also in
 {\sc R.\ Engui\c{c}a}, {\sc A.\ Gavioli}, and {\sc L.\ Sanchez}
 \cite[Theorem 4.2, p.~182]{EGSanchez}
with a sketched proof.

%%%%%%%%%%%%%%%%%%%%%%%%%%%%%%%%%%%%%%%%%%%%%%%%%%%%%%%%%%%%%%%%%%%%%%%%%%%
\section{Uniqueness Result for $c^{\ast}$}
\label{s:Unique}

For the sake of reader's convenience,
in this section we establish a few comparison results
(weak and strong)
which provide important technical tools.
These results are essentially due to {\sc E.\ Kamke},
as we have already mentioned in the proof of
Proposition~\ref{prop-c^*} above.
They will enable us to establish the uniqueness of $c^{\ast}$
specified in Corollary~\ref{cor-c<c^*}
(see Corollary~\ref{cor-uniq_c^*} below).

%%%%%%%%%%%%%%%%%%%%%%%%%%%%%%%%%%%%%%%%%%%%%%%%%%%%%%%%%%%%%%%%%%%%%%%
%%%%%    Sub- and Supersolution (Definition)    %%%%%%%%%%%%%%%%%%%%%%%
%%%%%%%%%%%%%%%%%%%%%%%%%%%%%%%%%%%%%%%%%%%%%%%%%%%%%%%%%%%%%%%%%%%%%%%
\begin{definition}\label{def-sub/super}\nopagebreak
\begingroup\rm
Let us consider the differential equation \eqref{eq:FKPP:y(r)}
on some nonempty open interval
$(a,b)\subset (-1,1)\subset \RR$, regardless
an {\em initial condition\/} at $a$ or
a {\em terminal condition\/} at~$b$.

We say that a continuously differentiable function
$\underline{y}\colon [a,b]\to \RR$
is a {\em subsolution\/} to equation~\eqref{eq:FKPP:y(r)}
on $(a,b)$ if and only if
\begin{equation}
\label{sub:FKPP:y(-1)}
    \frac{ \mathrm{d} \underline{y} }{\mathrm{d}r}
  \leq p'\left( c\, ( \underline{y}^{+} )^{1/p} + g(r) \right) \,,
    \quad r\in (a,b) \,.
\end{equation}
Similarly, we say that a continuously differentiable function
$\overline{y}\colon [a,b]\to \RR$
is a {\em supersolution\/} to equation~\eqref{eq:FKPP:y(r)}
on $(a,b)$ if and only if
\begin{equation}
\label{sup:FKPP:y(-1)}
    \frac{ \mathrm{d} \overline{y} }{\mathrm{d}r}
  \geq p'\left( c\, ( \overline{y}^{+} )^{1/p} + g(r) \right) \,,
    \quad r\in (a,b) \,.
\end{equation}
\endgroup
\end{definition}
%%%%%%%%%%%%%%%%%%%%%%%%%%%%%%%%%%%%%%%%%%%%%%%%%%%%%%%%%%%%%%%%%%%%%%%
\par\vskip 10pt

%%%%%%%%%%%%%%%%%%%%%%%%%%%%%%%%%%%%%%%%%%%%%%%%%%%%%%%%%%%%%%%%%%%%%%%
%%%%%    Sub- and Supersolution (Remark)    %%%%%%%%%%%%%%%%%%%%%%%%%%%
%%%%%%%%%%%%%%%%%%%%%%%%%%%%%%%%%%%%%%%%%%%%%%%%%%%%%%%%%%%%%%%%%%%%%%%
\begin{remark}\label{rem-sub/super}\nopagebreak
\begingroup\rm
We remark that the terminology for sub- and super\-solutions
introduced above fits the case when eq.~\eqref{eq:FKPP:y(r)}
is supplemented by an initial condition for $y(a)$
at the initial point $r=a$.
However, if eq.~\eqref{eq:FKPP:y(r)}
is supplemented by a terminal condition for $y(b)$
at the terminal point $r=b$, it would be logically more appropriate
to switch the inequalities
\eqref{sub:FKPP:y(-1)} and \eqref{sup:FKPP:y(-1)}
in the definitions of sub- and super\-solutions; compare
{\rm Parts (a)} and~{\rm (b)} of Proposition~\ref{prop-y:sub/sup} below.
Namely, a natural way to view a terminal value problem is to revert
the (``time'') variable $r\mapsto -r$ which entails also
the reversion of inequalities
\eqref{sub:FKPP:y(-1)} and \eqref{sup:FKPP:y(-1)}.
However, we will not perform this reversion to avoid possible confusion.
\endgroup
\end{remark}
%%%%%%%%%%%%%%%%%%%%%%%%%%%%%%%%%%%%%%%%%%%%%%%%%%%%%%%%%%%%%%%%%%%%%%%
%\par\vskip 10pt

We begin with the following comparison result.

%%%%%%%%%%%%%%%%%%%%%%%%%%%%%%%%%%%%%%%%%%%%%%%%%%%%%%%%%%%%%%%%%%%%%%%
%%%%%    Sub- and Supersolution (Proposition)    %%%%%%%%%%%%%%%%%%%%%%
%%%%%%%%%%%%%%%%%%%%%%%%%%%%%%%%%%%%%%%%%%%%%%%%%%%%%%%%%%%%%%%%%%%%%%%
\begin{proposition}\label{prop-y:sub/sup}
Let\/ $c\in \RR$, $-1\leq a < b\leq 1$, and assume that\/
$\underline{y}, \overline{y}\colon [a,b]\to \RR$
is a pair of continuously differentiable functions, such that\/
$\underline{y}$ ($\overline{y}$)
is a subsolution $($supersolution, respectively$)$
to equation~\eqref{eq:FKPP:y(r)} on $(a,b)$.
\begin{itemize}
\item[{\rm (a)}]
If\/ $c\leq 0$ and\/
$\underline{y}(a)\leq \overline{y}(a)$ then we have also
$\underline{y}(r)\leq \overline{y}(r)$ for all\/ $r\in [a,b]$.
\item[{\rm (b)}]
If\/ $c\geq 0$ and\/
$\underline{y}(b)\geq \overline{y}(b)$ then we have also
$\underline{y}(r)\geq \overline{y}(r)$ for all\/ $r\in [a,b]$.
\end{itemize}
\end{proposition}
%%%%%%%%%%%%%%%%%%%%%%%%%%%%%%%%%%%%%%%%%%%%%%%%%%%%%%%%%%%%%%%%%%%%%%%
\par\vskip 10pt

Recall our abbreviations
$\RR_+\eqdef [0,\infty)$ and $\RR_-\eqdef (-\infty,0]$.

\par\vskip 10pt
%%%%%%%%%%%%%%%%%%%%%%%%%%%%
\proof
We prove only {\rm Part~(a)};
{\rm Part~(b)} is proved analogously.

Hence, assume that $c\leq 0$.
We subtract ineq.~\eqref{sup:FKPP:y(-1)} from \eqref{sub:FKPP:y(-1)},
thus arriving at
\begin{equation*}
%\label{sub-sup:FKPP:y(-1)}
    \frac{\mathrm{d}}{\mathrm{d}r}\, (\underline{y} - \overline{y})
  \leq p'c
    \left( (\underline{y}^{+})^{1/p} - (\overline{y}^{+})^{1/p} \right) \,,
    \quad r\in (a,b) \,.
\end{equation*}
Now let us multiply this inequality by
$(\underline{y} - \overline{y})^{+}$, thus obtaining
\begin{equation}
\label{sub-sup:y(-1)}
    (\underline{y} - \overline{y})^{+}\,
    \frac{\mathrm{d}}{\mathrm{d}r}\, (\underline{y} - \overline{y})
  \leq p'c
    \left( (\underline{y}^{+})^{1/p} - (\overline{y}^{+})^{1/p} \right)
    (\underline{y} - \overline{y})^{+}
  \leq 0 \,,\quad r\in (a,b) \,,
\end{equation}
thanks to $c\leq 0$ combined with the montonicity of the functions
$\xi\mapsto \xi^{+}$ and $\xi\mapsto (\xi^{+})^{1/p}$
from $\RR$ to $\RR_+$.
The last inequality entails
\begin{equation*}
%\label{sub-sup:y(-1)}
  \frac{1}{2}\,
    \frac{\mathrm{d}}{\mathrm{d}r}
    \left( (\underline{y} - \overline{y})^{+} \right)^2
  \leq 0 \quad\mbox{ for almost every }\, r\in (a,b) \,,
\end{equation*}
which shows that
$w\eqdef \left( (\underline{y} - \overline{y})^{+} \right)^2$
is a monotone nonincreasing, nonnegative function on $[a,b]$.
Since $w(a) = 0$ by our hypothesis, it follows that
$w(r) = 0$ holds for every $r\in [a,b]$.
This completes our proof of {\rm Part~(a)}.
\qed
%%%%%%%%%%%%%%%%%%%%%%%%%%%%
\par\vskip 10pt

The following corollary on the uniqueness and
the monotone dependence on the parameter $c\in \RR_-$
of the solution $y\equiv y_c\colon [-1,1]\to \RR$
to the initial value problem for eq.~\eqref{eq:FKPP:y(r)}
with the initial condition $y(-1) = 0$
is an easy direct consequence of Proposition~\ref{prop-y:sub/sup}.

%%%%%%%%%%%%%%%%%%%%%%%%%%%%%%%%%%%%%%%%%%%%%%%%%%%%%%%%%%%%%%%%%%%%%%%
%%%%%    c --> y_c(u) is monotone increasing (Corollary)    %%%%%%%%%%%
%%%%%%%%%%%%%%%%%%%%%%%%%%%%%%%%%%%%%%%%%%%%%%%%%%%%%%%%%%%%%%%%%%%%%%%
\begin{corollary}\label{cor-y_c:monot}
{\rm (i)}$\;$
Given any\/ $c\in \RR_-$,
the initial value problem for eq.~\eqref{eq:FKPP:y(r)}
with the initial condition $y(-1) = 0$
possesses a unique $($continuously differentiable$)$ solution
$y\equiv y_c\colon [-1,1]\to \RR$.

{\rm (ii)}$\;$
If\/ $-\infty < c_1\leq c_2\leq 0$ then
$y_{c_1}\leq y_{c_2}$ holds pointwise throughout $[-1,1]$.
In particular, the function
$c\mapsto y_c(1)\colon \RR_-\to \RR$ is monotone increasing.
\end{corollary}
%%%%%%%%%%%%%%%%%%%%%%%%%%%%%%%%%%%%%%%%%%%%%%%%%%%%%%%%%%%%%%%%%%%%%%%
%\par\vskip 10pt

\par\vskip 10pt
%%%%%%%%%%%%%%%%%%%%%%%%%%%%
\proof
{\rm Part~(i)} (the uniqueness claim) follows trivially from
Proposition~\ref{prop-y:sub/sup}, {\rm Part~(a)}.

{\rm Part~(ii)}:
Let $-\infty < c_1\leq c_2\leq 0$ and choose any constant\/
$c\in [c_1,c_2]$; hence, $c\leq 0$.
Then $y_{c_1}, y_{c_2}\colon [-1,1]\to \RR$
is a pair of sub- and super\-solutions
to equation~\eqref{eq:FKPP:y(r)} on $(-1,1)$, respectively, with
$c$ chosen above ($c_1\leq c\leq c_2$)
and the initial condition
$y_{c_1}(-1) = y_{c_2}(-1) = 0$.
The pointwise ordering
$y_{c_1}\leq y_{c_2}$ throughout $[-1,1]$ now follows from
Proposition~\ref{prop-y:sub/sup}, {\rm Part~(a)}, again.
\qed
%%%%%%%%%%%%%%%%%%%%%%%%%%%%
\par\vskip 10pt

We complement Proposition~\ref{prop-y:sub/sup}
with another comparison result.

%%%%%%%%%%%%%%%%%%%%%%%%%%%%%%%%%%%%%%%%%%%%%%%%%%%%%%%%%%%%%%%%%%%%%%%
%%%%%    Sub- and Supersolution (Proposition)    %%%%%%%%%%%%%%%%%%%%%%
%%%%%%%%%%%%%%%%%%%%%%%%%%%%%%%%%%%%%%%%%%%%%%%%%%%%%%%%%%%%%%%%%%%%%%%
\begin{proposition}\label{prop-y(1):sub/sup}
Let\/ $c\in \RR$, $a\in [s_0,1)$, and assume that\/
$\underline{y}, \overline{y}\colon [a,1]\to \RR$
is a pair of continuously differentiable functions, such that\/
$\underline{y}$ ($\overline{y}$)
is a subsolution $($supersolution, respectively$)$
to equation~\eqref{eq:FKPP:y(r)} on $(a,1)$.
Moreover, assume that\/
\begin{itemize}
\item[{}]
$\underline{y}(r) > 0$ and\/ $\overline{y}(r) > 0$ for all\/ 
$r\in [a,1)$, together with\/
$\overline{y}(1) = 0\leq \underline{y}(1)$.
\end{itemize}
Then we have also
\begin{itemize}
\item[{\rm (a)}]
$\overline{y}(r)\leq \underline{y}(r)$ for all\/ $r\in [a,1]$.
\item[{\rm (b)}]
More precisely, if\/
$\overline{y}(1) = 0 < \underline{y}(1)$ then
$\overline{y}(r) < \underline{y}(r)$ holds for all\/ $r\in [a,1]$. 
If\/
$\overline{y}(1) = \underline{y}(1) = 0$ then there is some number\/
$a'\in [a,1]$ such that\/
$\overline{y}(r) < \underline{y}(r)$ for all\/ $r\in [a,a')$ and\/
$\overline{y}(r) = \underline{y}(r)$ for all\/ $r\in [a',1]$.
\end{itemize}
\end{proposition}
%%%%%%%%%%%%%%%%%%%%%%%%%%%%%%%%%%%%%%%%%%%%%%%%%%%%%%%%%%%%%%%%%%%%%%%
%\par\vskip 10pt

The comparison result in Part~{\rm (b)}
takes care of the lack of uniqueness in
the terminal value problem for eq.~\eqref{eq:FKPP:y(r)}
with the terminal condition $y(1) = 0$ and $c<0$, as announced in
Remark~\ref{rem-sol:G(u)}.

\par\vskip 10pt
%%%%%%%%%%%%%%%%%%%%%%%%%%%%
{\it Proof of\/} Proposition~\ref{prop-y(1):sub/sup}.
Part~{\rm (a)}:
We begin by substituting the variable
$t = 1-r$ with $0\leq t\leq 1-a$ $(\leq 1 - s_0)$ and the functions
\begin{equation*}
  \underline{z}(t)\eqdef (\underline{y}(1-t))^{1/p'}
    \qquad\mbox{ and }\qquad
  \overline{z}(t) \eqdef (\overline{y}(1-t))^{1/p'}
\end{equation*}
with the derivatives
(recall that $\frac{1}{p} + \frac{1}{p'} = 1$)
\begin{align*}
  \frac{\mathrm{d}}{\mathrm{d}t}\, \underline{z}(t)
& = {}- \frac{1}{p'}\, (\underline{y}(1-t))^{-1/p}\cdot
  \frac{\mathrm{d}}{\mathrm{d}t}\, \underline{y}(1-t)
    \quad\mbox{ and }\quad
\\
  \frac{\mathrm{d}}{\mathrm{d}t}\, \overline{z}(t)
& = {}- \frac{1}{p'}\, (\overline{y}(1-t))^{-1/p}\cdot
  \frac{\mathrm{d}}{\mathrm{d}t}\, \overline{y}(1-t)
    \quad\mbox{ for }\, 0 < t\leq 1-a \,.
\end{align*}
Hence, inequalities
\eqref{sub:FKPP:y(-1)} and \eqref{sup:FKPP:y(-1)}, respectively,
are equivalent with
\begin{equation}
\label{sub:FKPP:z(1)}
    \frac{ \mathrm{d} \underline{z} }{\mathrm{d}t}
  \geq {}- c - g(1-t)\, ( \underline{z}(t) )^{-1/(p-1)} \,,
    \quad t\in (0,1-a] \,,
\end{equation}
and
\begin{equation}
\label{sup:FKPP:z(1)}
    \frac{ \mathrm{d} \overline{z} }{\mathrm{d}t}
  \leq {}- c - g(1-t)\, ( \overline{z}(t) )^{-1/(p-1)} \,,
    \quad t\in (0,1-a] \,.
\end{equation}

We subtract ineq.~\eqref{sub:FKPP:z(1)} from \eqref{sup:FKPP:z(1)},
thus arriving at
\begin{equation}
\label{sub-sup:FKPP:z(1)}
    \frac{\mathrm{d}}{\mathrm{d}t}\, (\overline{z} - \underline{z})
  \leq {}- g(1-t)
    \left[ ( \overline{z}(t) )^{-1/(p-1)}
         - ( \underline{z}(t) )^{-1/(p-1)}
    \right] \,,
    \quad t\in (0,1-a] \,.
\end{equation}
Now let us multiply ineq.\ \eqref{sub-sup:FKPP:z(1)} by
$(\overline{z} - \underline{z})^{+}$, thus obtaining
\begin{equation}
\label{sub-sup:z(1)}
    (\overline{z} - \underline{z})^{+}\,
    \frac{\mathrm{d}}{\mathrm{d}t}\, (\overline{z} - \underline{z})
  \leq {}- g(1-t)
    \left[ ( \overline{z}(t) )^{-1/(p-1)}
         - ( \underline{z}(t) )^{-1/(p-1)}
    \right]
    (\overline{z} - \underline{z})^{+}
  \leq 0
\end{equation}
for all $t\in (0,1-a]$, thanks to
$g(1-t)\leq 0$ for $0 < t\leq 1-a\leq 1 - s_0$.
The last inequality entails
\begin{equation*}
%\label{sub-sup:y(-1)}
  \frac{1}{2}\,
    \frac{\mathrm{d}}{\mathrm{d}t}
    \left( (\overline{z} - \underline{z})^{+} \right)^2
  \leq 0 \quad\mbox{ for almost every }\, r\in (0,1-a) \,,
\end{equation*}
which shows that
$w\eqdef \left( (\overline{z} - \underline{z})^{+} \right)^2$
is a monotone nonincreasing, nonnegative function on $[0,1-a]$.
Since $w(0) = 0$ by our hypothesis, it follows that
$w(t) = 0$ holds for every $t\in [0,1-a]$.
Equivalently,
$\overline{y}(r)\leq \underline{y}(r)$ holds for all $r\in [a,1]$,
which proves Part~{\rm (a)}.

Part~{\rm (b)}:
We return to ineq.~\eqref{sub-sup:FKPP:z(1)}.
The difference on the right\--hand side takes the form
\begin{equation}
\label{est:z^-z_}
\begin{aligned}
& 0\leq
    ( \overline{z}(t) )^{-1/(p-1)}
  - ( \underline{z}(t) )^{-1/(p-1)}
  =
\\
& {}- \frac{1}{p-1}
    \left( \int_0^1
    \left[ (1-\theta) \overline{z}(t) + \theta \underline{z}(t)
    \right]^{-p'} \,\mathrm{d}\theta
    \right)
    \left( \overline{z}(t) - \underline{z}(t) \right)
\\
& \leq L\left( \underline{z}(t) - \overline{z}(t) \right)
  \quad\mbox{ for all }\, t\in [1-b,1-a] \,,
\end{aligned}
\end{equation}
where $b\in (a,1)$ is arbitrary and
the Lipschitz constant $L\equiv L(b)\in (0,\infty)$ depends on $b$.
We may take
\begin{equation*}
  L = \frac{1}{p-1}
      \left( \min_{ t\in [1-b,1-a] } \overline{z}(t) \right)^{-p'}
    \in (0,\infty) \,,
\end{equation*}
by the inequality
\begin{equation*}
  \int_0^1
    \left[ (1-\theta) \overline{z}(t) + \theta \underline{z}(t)
    \right]^{-p'} \,\mathrm{d}\theta
  \leq ( \overline{z}(t) )^{-p'} \,,
    \quad t\in (0,1-a] \,.
\end{equation*}
Consequently, we can estimate the right\--hand side of
ineq.~\eqref{sub-sup:FKPP:z(1)} as follows:
\begin{equation}
\label{sub-sup:z^-z_}
    \frac{\mathrm{d}}{\mathrm{d}t}\, (\overline{z} - \underline{z})
  \leq \hat{L}
    \left( \underline{z}(t) - \overline{z}(t) \right) \,,
    \quad t\in [1-b,1-a] \,,
\end{equation}
where
\begin{equation*}
  \hat{L}\equiv \hat{L}(b)\eqdef L\cdot \sup_{[-1,1]} |g| \,,\qquad
  \hat{L}\in (0,\infty) \,.
\end{equation*}
Ineq.~\eqref{sub-sup:z^-z_} is equivalent with
\begin{equation}
\label{exp:z^-z_}
    \frac{\mathrm{d}}{\mathrm{d}t}
    \left[ \mathrm{e}^{\hat{L} t}
    \left( \overline{z}(t) - \underline{z}(t) \right)
    \right] \leq 0
    \quad\mbox{ for every }\, t\in [1-b,1-a] \,.
\end{equation}
By integration over any compact interval
$[t_1,t_2]\subset [1-b,1-a]$, we get
\begin{equation}
\label{a<b:z^-z_}
    \overline{z}(t_2) - \underline{z}(t_2)
  \leq \mathrm{e}^{- \hat{L} (t_2 - t_1)}
    \left( \overline{z}(t_1) - \underline{z}(t_1) \right)
  \leq 0 \,.
\end{equation}

First, assume that
$\overline{y}(1) = 0 < \underline{y}(1)$, that is,
$\overline{z}(0) = 0 < \underline{z}(0)$.
Then, by continuity, also
$\overline{z}(t_1) < \underline{z}(t_1)$
holds provided $t_1\in (0,1-a)$ is chosen to be small enough.
From ineq.~\eqref{a<b:z^-z_} we deduce that also
$\overline{z}(t) < \underline{z}(t)$ for all $t\in [t_1,1-a]$.
We have proved
$\overline{z}(t) < \underline{z}(t)$ for every $t\in [0,1-a]$.
The desired inequality
$\overline{y}(r) < \underline{y}(r)$ for every $r\in [a,1]$
follows immediately.

Now assume
$\overline{y}(1) = \underline{y}(1) = 0$, that is,
$\overline{z}(0) = \underline{z}(0) = 0$.
Let $t'\in [0,1-a]$ be the greatest number $t\in [0,1-a]$, such that
$\overline{z}(t) - \underline{z}(t) = 0$
holds for every $t\in [0,t']$.
Equivalently,
$\overline{y}(r) - \underline{y}(r) = 0$
holds for every $r\in [a',1]$, where $a'= 1 - t'\in [a,1]$.
If $t'= 1-a$ then we are done.
Thus, assume $0\leq t'< 1-a$.
Then there is a sequence
$\{ \tau_n\}_{n=1}^{\infty} \subset (t',1-a)$, such that
$\tau_n\searrow t'$ as $n\nearrow \infty$ and
$\overline{z}(\tau_n) - \underline{z}(\tau_n) < 0$
for every $n = 1,2,3,\dots$.
Again, from ineq.~\eqref{a<b:z^-z_} we deduce that also
$\overline{z}(t) < \underline{z}(t)$ for all $t\in [\tau_n,1-a]$.
We have proved
$\overline{z}(t) < \underline{z}(t)$ for every $t\in (t',1-a]$.
The desired inequality
$\overline{y}(r) < \underline{y}(r)$ for every $r\in [a,a')$
follows immediately.

The proof is complete.
\qed
%%%%%%%%%%%%%%%%%%%%%%%%%%%%
\par\vskip 10pt

Corollary~\ref{cor-c<c^*}, Corollary~\ref{cor-y_c:monot}, and
Proposition~\ref{prop-y(1):sub/sup}
entail the following important result.

%%%%%%%%%%%%%%%%%%%%%%%%%%%%%%%%%%%%%%%%%%%%%%%%%%%%%%%%%%%%%%%%%%%%%%%
%%%%%    y(u) >< G(u) is a solution (Corollary)    %%%%%%%%%%%%%%%%%%%%
%%%%%%%%%%%%%%%%%%%%%%%%%%%%%%%%%%%%%%%%%%%%%%%%%%%%%%%%%%%%%%%%%%%%%%%
\begin{corollary}\label{cor-uniq_c^*}
There is a unique number\/ $c\in \RR$ such that
problem~\eqref{BVP:FKPP:y(r)} possesses a solution
$y\colon [-1,1]\to \RR$, namely, $c = c^{\ast}$ $(< 0)$
defined in eq.~\eqref{def:c^*}.
This solution is given by $y = y_{c^{\ast}}$ and satisfies
$y(r) > 0$ for all\/ $r\in (-1,1)$.
\end{corollary}
%%%%%%%%%%%%%%%%%%%%%%%%%%%%%%%%%%%%%%%%%%%%%%%%%%%%%%%%%%%%%%%%%%%%%%%
%\par\vskip 10pt

\par\vskip 10pt
%%%%%%%%%%%%%%%%%%%%%%%%%%%%
\proof
The existence of $c^{\ast}$ $(< 0)$ has been established in
Corollary~\ref{cor-c<c^*}.
On the contrary to uniqueness, suppose that there is another number
$c$, say, $c = c'\in \RR$, $c'\neq c^{\ast}$, such that
$y = y_{c'}$ satisfies the boundary conditions
$y(-1) = y(1) = 0$.
We have $y_{c'}(r) > 0$ and $y_{c^{\ast}}(r) > 0$ for all $r\in (-1,1)$,
by Lemma~\ref{lem-y_c>0}, Part~{\rm (ii)}.
By our definition of $c^{\ast}$ in eq.~\eqref{def:c^*},
we must have $c'< c^{\ast}$ $(< 0)$.
Fixing any constant $c\in [c', c^{\ast}]$, we observe that
$y_{c'}$ ($y_{c^{\ast}}$, respectively)
is a subsolution (supersolution)
to equation~\eqref{eq:FKPP:y(r)} on $(-1,1)$ (with this $c$).

Corollary~\ref{cor-y_c:monot}, Part~{\rm (ii)}, implies
$y_{c'}\leq y_{c^{\ast}}$ throughout $[-1,1]$.
In contrast, Proposition~\ref{prop-y(1):sub/sup}, Part~{\rm (b)},
forces $y_{c^{\ast}}\leq y_{c'}$ on $[s_0,1]$.
We conclude that
$y_{c'}\equiv y_{c^{\ast}}$ on $[s_0,1]$.
Consequently, both,
$y_{c'}$ and $y_{c^{\ast}}$, verify eq.~\eqref{eq:FKPP:y(r)}
on the interval $(s_0,1)$ for both values of $c$,
$c = c'$ and $c = c^{\ast}$.
Since also $y_{c'}(r) = y_{c^{\ast}}(r) > 0$ for every $r\in [s_0,1)$,
eq.~\eqref{eq:FKPP:y(r)}
with $c = c'$ and $c = c^{\ast}$ forces $c'= c^{\ast}$.
We have proved also the uniqueness.
\qed
%%%%%%%%%%%%%%%%%%%%%%%%%%%%
\par\vskip 10pt

%%%%%%%%%%%%%%%%%%%%%%%%%%%%%%%%%%%%%%%%%%%%%%%%%%%%%%%%%%%%%%%%%%%%%%%%%%%
\section{Asymptotic Behavior}
\label{s:Asympt}

In this section we investigate
the {\it asymptotic behavior\/} of the (unique) solution
$y = y_{c^{\ast}}\colon [-1,1]\to \RR$
to problem~\eqref{BVP:FKPP:y(r)} with $c = c^{\ast}$.
Of course, we assume $G(1) > 0$ throughout the present section.
We recall that the existence and uniqueness of both,
$c^{\ast}$ $(< 0)$ and $y_{c^{\ast}}$, have been obtained in
Corollaries \ref{cor-c<c^*} and~\ref{cor-uniq_c^*}.

In fact, the asymptotic estimates for $r\in (-1,1)$ near the endpoints
$\pm 1$ obtained in this section remain valid in
the following more general situation with an arbitary value of 
the parameter $c\leq 0$ near the left endpoint $-1$:
We assume that
$y = y_c\colon [-1,1]\to \RR$ is the (unique) solution
to eq.~\eqref{eq:FKPP:y(r)}
satisfying the initial condition $y(-1) = 0$.

We assume that
$g(r)\eqdef d(r)^{1/(p-1)}\, f(r)$ satisfies also
the following hypothesis:
\begin{itemize}
\item[{}]
There are constants $\gamma, \gamma_0\in (0,\infty)$
such that
\begin{equation}
\label{lim:g(r)}
  \lim_{r\to -1}\, \frac{g(r)}{ (1+r)^{\gamma} } = \gamma_0 \,.
\end{equation}
\end{itemize}
Equivalently, we have
\begin{equation}
%\label{lim:g(r)}
\nonumber
\tag{\ref{lim:g(r)}$'$}
  g(r) = ( \gamma_0 + \eta(r) ) (1+r)^{\gamma}
  \quad\mbox{ for }\, r\in [-1,1] \,,
\end{equation}
where $\eta\colon [-1,1]\to \RR$ is a continuous function with
$\eta(0) = 0$.

We abbreviate the differential operator
\begin{equation}
\label{def:A(r)}
  (\mathcal{A}y)(r)\eqdef
    \frac{\mathrm{d}y}{\mathrm{d}r}
  - p'\left( c\, (y^{+})^{1/p} + g(r) \right) \,,
    \quad r\in (-1,1) \,,
\end{equation}
for any function $y\in C^1([-1,1])$.
We are interested in constructing sub- and supersolutions
to equation~\eqref{eq:FKPP:y(r)} as defined in Section~\ref{s:Unique}.
Let $\kappa\in (0,\infty)$ be a constant (to be determined later)
and let us consider the (nonnegative) function
\begin{equation}
\label{e:w_kappa}
  w_{\kappa}(r)\eqdef \kappa\, (1+r)^{1+\gamma}
  \quad\mbox{ for }\, r\in [-1,1] \,.
\end{equation}
Then we have
\begin{equation}
\label{e:Aw_kappa}
\begin{aligned}
& (\mathcal{A}w_{\kappa})(r)
  = \frac{ \mathrm{d}w_{\kappa} }{\mathrm{d}r}
  - p'\left( c\, (w_{\kappa})^{1/p} + g(r) \right)
\\
& = \kappa (1 + \gamma) (1+r)^{\gamma}
  - p'\left( c\, \kappa^{1/p} (1+r)^{(1+\gamma)/p}
           + ( \gamma_0 + \eta(r) ) (1+r)^{\gamma}
      \right)
\\
& = \left[ \kappa (1 + \gamma) - p'\, ( \gamma_0 + \eta(r) )
    \right] (1+r)^{\gamma}
  - p' c\, \kappa^{1/p} (1+r)^{(1+\gamma)/p}
\end{aligned}
\end{equation}
for $r\in (-1,1)$.
Letting $r\to -1$ we observe that the first expression, with the power
$(1+r)^{\gamma}$, dominates the second one, with the power 
$(1+r)^{(1+\gamma)/p}$, if and only if
$\gamma\leq (1+\gamma)/p$ holds, that is,
$\gamma\leq p'/p = 1/(p-1)$.
Consequently, if this is the case, then,
for every $r\in (-1,1)$ sufficiently close to $-1$ we have:
\begin{itemize}
\item[{}]
$(\mathcal{A}w_{\kappa})(r) < 0$ provided $\kappa > 0$ is small enough,
and
\item[{}]
$(\mathcal{A}w_{\kappa})(r) > 0$ provided $\kappa > 0$ is large enough.
\end{itemize}
We can state these simple facts about
$w_{\kappa}$ being a sub- or supersolution to problem \eqref{eq:FKPP:y(r)}
on some interval $(-1,-1+\varrho)$ as follows ($0 < \varrho < 2$):

%%%%%%%%%%%%%%%%%%%%%%%%%%%%%%%%%%%%%%%%%%%%%%%%%%%%%%%%%%%%%%%%%%%%%%%
%%%%%    y(u) >< G(u) is a solution (Lemma)    %%%%%%%%%%%%%%%%%%%%%%%%
%%%%%%%%%%%%%%%%%%%%%%%%%%%%%%%%%%%%%%%%%%%%%%%%%%%%%%%%%%%%%%%%%%%%%%%
\begin{lemma}\label{lem-power<g(r)}
Let\/ $c\in \RR$ and assume that\/
$g\colon [-1,1]\to \RR$ satisfies \eqref{lim:g(r)}
with $0 < \gamma\leq 1/(p-1)$.
Then there exist numbers $\varrho\in (0,2)$ and\/
$0 < \underline{\kappa} < \overline{\kappa} < \infty$,
such that
\begin{enumerate}
\renewcommand{\labelenumi}{(\roman{enumi})}
\item[{\rm (i)}]
if\/ $0 < \kappa\leq \underline{\kappa}$ then
$(\mathcal{A}w_{\kappa})(r) < 0$ holds for all\/
$r\in (-1,-1+\varrho)$; and
\item[{\rm (ii)}]
if\/ $\overline{\kappa}\leq \kappa < \infty$ then
$(\mathcal{A}w_{\kappa})(r) > 0$ holds for all\/
$r\in (-1,-1+\varrho)$.
\end{enumerate}
\end{lemma}
%%%%%%%%%%%%%%%%%%%%%%%%%%%%%%%%%%%%%%%%%%%%%%%%%%%%%%%%%%%%%%%%%%%%%%%
\par\vskip 10pt

We remark that all three numbers
$\varrho$, $\underline{\kappa}$, and $\overline{\kappa}$
may depend on $c$, in general.

Now, let us fix $c = c^{\ast}$ given by eq.~\eqref{def:c^*}
and recall Corollary~\ref{cor-uniq_c^*}.
We combine Proposition~\ref{prop-y:sub/sup}, Part~{\rm (a)},
with Lemma~\ref{lem-power<g(r)}
to conclude that
\begin{equation}
\label{ineq:y_c*}
  w_{\underline{\kappa}}(r) = \underline{\kappa}\, (1+r)^{1+\gamma}
  \leq y_{c^{\ast}}(r)\leq
  w_{\overline{\kappa}}(r) = \overline{\kappa}\, (1+r)^{1+\gamma}
\end{equation}
for all $r\in (-1,-1+\varrho)$.
We recall that
$V(U) = \left[ y_{c^{\ast}}(r)\right]^{1/p'}$
with $U=r$, by eq.~\eqref{e:y=V^p'}.
Consequently, inequalities \eqref{ineq:y_c*} above read
\begin{equation}
\label{ineq:V_c*}
  (\underline{\kappa})^{1/p'}\, (1+U)^{(1+\gamma) / p'}
  \leq V(U)\leq
  (\overline{\kappa})^{1/p'}\, (1+U)^{(1+\gamma) / p'}
\end{equation}
for all $U\in (-1,-1+\varrho)$.
In eq.~\eqref{e:dU/dx}, that is
\begin{equation}
%\label{e:dx/dU}
\nonumber
\tag{\ref{e:dx/dU}}
    \genfrac{}{}{}0{\mathrm{d}x}{\mathrm{d}U}
  = {}-
    \genfrac{(}{)}{}0{d(U)}{V(U)}^{1/(p-1)} < 0 \,,\quad
    -1 < U < 1 \,,
\end{equation}
for the unknown function $x = x(U)$, we need
the following equivalent form of \eqref{ineq:V_c*},
\begin{equation}
%\label{ineq:V_c*}
\nonumber
\tag{\ref{ineq:V_c*}$'$}
  (\underline{\kappa})^{1/p}\, (1+U)^{(1+\gamma) / p}
  \leq V(U)^{1/(p-1)}\leq
  (\overline{\kappa})^{1/p}\, (1+U)^{(1+\gamma) / p}
\end{equation}
for all $U\in (-1,-1+\varrho)$.

It remains to investigate the case $\gamma > (1+\gamma)/p$, that is,
$\gamma > p'/p = 1/(p-1)$, in which the second expression
in eq.~\eqref{e:Aw_kappa}, with the power
$(1+r)^{(1+\gamma)/p}$, dominates the first one, with the power
$(1+r)^{\gamma}$.
Here, we assume $c < 0$.
Then, given any $\kappa\in (0,\infty)$, we have
(thanks to $c < 0$)
\begin{itemize}
\item[{}]
$(\mathcal{A}w_{\kappa})(r) > 0$
for every $r\in (-1,1)$ sufficiently close to $-1$.
\end{itemize}
We can state this simple fact about
$w_{\kappa}$ being a supersolution to problem \eqref{eq:FKPP:y(r)}
on some interval $(-1,-1+\varrho)$ as follows ($0 < \varrho < 2$):

%%%%%%%%%%%%%%%%%%%%%%%%%%%%%%%%%%%%%%%%%%%%%%%%%%%%%%%%%%%%%%%%%%%%%%%
%%%%%    y(u) >< G(u) is a solution (Lemma)    %%%%%%%%%%%%%%%%%%%%%%%%
%%%%%%%%%%%%%%%%%%%%%%%%%%%%%%%%%%%%%%%%%%%%%%%%%%%%%%%%%%%%%%%%%%%%%%%
\begin{lemma}\label{lem-power>g(r)}
Let\/ $c < 0$ and assume that\/
$g\colon [-1,1]\to \RR$ satisfies \eqref{lim:g(r)}
with $\gamma > 1/(p-1)$.
Then, given any $\kappa\in (0,\infty)$,
there exists a number $\varrho\in (0,2)$, such that
\begin{itemize}
\item[{}]
$(\mathcal{A}w_{\kappa})(r) > 0$ holds for all\/
$r\in (-1,-1+\varrho)$.
\end{itemize}
\end{lemma}
%%%%%%%%%%%%%%%%%%%%%%%%%%%%%%%%%%%%%%%%%%%%%%%%%%%%%%%%%%%%%%%%%%%%%%%
\par\vskip 10pt

We remark that, again, the number $\varrho$
may depend on $\kappa$ and $c$, in general.

Now, let us fix $c = c^{\ast}$ given by eq.~\eqref{def:c^*}
and recall Corollary~\ref{cor-uniq_c^*}.
We combine Proposition~\ref{prop-y:sub/sup}, Part~{\rm (a)},
with Lemma~\ref{lem-power<g(r)}
to conclude that
\begin{equation}
\label{inequ:y_c*}
  y_{c^{\ast}}(r)\leq w_{\kappa}(r) = \kappa\, (1+r)^{1+\gamma}
\end{equation}
for all $r\in (-1,-1+\varrho)$.
We recall that
$V(U) = \left[ y_{c^{\ast}}(r)\right]^{1/p'}$
with $U=r$, by eq.~\eqref{e:y=V^p'}.
Consequently, ineq.~\eqref{inequ:y_c*} above reads
\begin{equation}
\label{inequ:V_c*}
  V(U)\leq \kappa^{1/p'}\, (1+U)^{(1+\gamma) / p'}
\end{equation}
for all $U\in (-1,-1+\varrho)$.
In eq.~\eqref{e:dx/dU}, that is
\begin{equation*}
%\label{e:dx/dU}
    \genfrac{}{}{}0{\mathrm{d}x}{\mathrm{d}U}
  = {}-
    \genfrac{(}{)}{}0{d(U)}{V(U)}^{1/(p-1)} < 0 \,,\quad
    -1 < U < 1 \,,
\end{equation*}
for the unknown function $x = x(U)$, we need
the following equivalent form of \eqref{inequ:V_c*},
\begin{equation}
%\label{ineq:V_c*}
\nonumber
\tag{\ref{inequ:V_c*}$'$}
  V(U)^{1/(p-1)}\leq \kappa^{1/p}\, (1+U)^{(1+\gamma) / p}
\end{equation}
for all $U\in (-1,-1+\varrho)$.

%%%%%%%%%%%%%%%%%%%%%%%%%%%%%%%%%%%%%%%%%%%%%%%%%%%%%%%%%%%%%%%%%%%%%%%
%%%%%    y(u) >< G(u) is a solution (Corollary)    %%%%%%%%%%%%%%%%%%%%
%%%%%%%%%%%%%%%%%%%%%%%%%%%%%%%%%%%%%%%%%%%%%%%%%%%%%%%%%%%%%%%%%%%%%%%
\begin{corollary}\label{cor-power<g(r)}
Let\/ $c = c^{\ast}$ $(< 0)$ and assume that\/
$g\colon [-1,1]\to \RR$ satisfies \eqref{lim:g(r)}.
Then we have the following conclusions for the limit\/
$x_{-1}\eqdef \lim_{U\to -1} x(U)$ $({}\leq +\infty)$
of the solution $x = x(U)$ to the differential equation
\eqref{e:dx/dU}:
\begin{enumerate}
\renewcommand{\labelenumi}{(\roman{enumi})}
\item[{\rm (i)}]
In case $1 < p\leq 2$, the limit\/ $x_{-1}$ is finite if and only if\/
$0 < \gamma < p-1$.
$($Equivalently, $x_{-1} = +\infty$ if and only if\/ $\gamma\geq p-1$.$)$
\item[{\rm (ii)}]
In case $2 < p < \infty$, the limit\/ $x_{-1} = +\infty$ if\/
$\gamma\geq p-1$.
\end{enumerate}

Consequently, we have $x_{-1} = +\infty$ whenever
$0 < p-1\leq \gamma < \infty$.
\end{corollary}
%%%%%%%%%%%%%%%%%%%%%%%%%%%%%%%%%%%%%%%%%%%%%%%%%%%%%%%%%%%%%%%%%%%%%%%
\par\vskip 10pt

The {\it proof\/} of Part~{\rm (i)}
(Part~{\rm (ii)}, respectively)
follows directly from eq.~\eqref{e:dx/dU} combined with
ineq.~\eqref{ineq:V_c*}$'$ (ineq.~\eqref{inequ:V_c*}$'$).

In the linear diffusion case, i.e., for $p=2$,
$x_{-1}$ is finite if $\gamma < 1$, whereas
$x_{-1} = +\infty$ if $\gamma\geq 1$.

%%%%%%%%%%%%%%%%%%%%%%%%%%%%%%%%%%%%%%%%%%%%%%%%%%%%%%%%%%%%%%%%%%%%%%%
%%%%%    y(u) >< G(u) is a solution (Remark)    %%%%%%%%%%%%%%%%%%%%%%%
%%%%%%%%%%%%%%%%%%%%%%%%%%%%%%%%%%%%%%%%%%%%%%%%%%%%%%%%%%%%%%%%%%%%%%%
\begin{remark}\label{rem-power<g(r)}\nopagebreak
\begingroup\rm
Under a condition similar to \eqref{lim:g(r)},
but imposed at the right endpoint as $r\to +1$,
we obtain conclusions analogous to those in\/
Corollary~\ref{cor-power<g(r)} for the limit\/
$x_1\eqdef \lim_{U\to +1} x(U)$ $({}\geq -\infty)$.
\endgroup
\end{remark}
%%%%%%%%%%%%%%%%%%%%%%%%%%%%%%%%%%%%%%%%%%%%%%%%%%%%%%%%%%%%%%%%%%%%%%%
\par\vskip 10pt

%\vfill
%%%%%%%%%%%%%%%%%%%%%%%%%%%%%%%%%%%%%%%%%%%%%%%%%%%%%%%%%%%%%%%%%%%%%%%
%%%%%    ACKNOWLEDGMENTS    %%%%%%%%%%%%%%%%%%%%%%%%%%%%%%%%%%%%%%%%%%%
%%%%%%%%%%%%%%%%%%%%%%%%%%%%%%%%%%%%%%%%%%%%%%%%%%%%%%%%%%%%%%%%%%%%%%%
\subsection*{Acknowledgments}
\begin{small}
The work of Pavel Dr\'abek was supported in part by
the Grant Agency of the Czech Republic (GA\v{C}R)
under Grant {\#}$13-00863$S, and
the work of Peter Tak\'a\v{c} by
a grant from Deutsche Forschungs\-gemeinschaft (DFG, Germany)
under Grant {\#} TA~213/16--1.
\end{small}
%%%%%%%%%%%%%%%%%%%%%%%%%%%%%%%%%%%%%%%%%%%%%%%%%%%%%%%%%%%%%%%%%%%%%%%

%%%%%%%%%%%%%%%%%%%%%%%%%%%%%%%%%%%%%%%%%%%%%%%%%%%%%%%%%%%%%%%%%%%%%%%
%%%%%    BIBLIOGRAPHY    %%%%%%%%%%%%%%%%%%%%%%%%%%%%%%%%%%%%%%%%%%%%%%
%%%%%%%%%%%%%%%%%%%%%%%%%%%%%%%%%%%%%%%%%%%%%%%%%%%%%%%%%%%%%%%%%%%%%%%

%%%%%%%%%%%%%%%%%%%%%%%%%%%%%%%%%%%%%%%%%%%%%%%%%%%%%%%%%%%%%%%%%%%%%%%%%%%%
%
%\bibliographystyle{amsplain}        %% if within the text
%
\makeatletter \renewcommand{\@biblabel}[1]{\hfill#1.} \makeatother
%
%\ifx\undefined\bysame
%\newcommand{\bysame}{\leavevmode\hbox to3em{\hrulefill}\,}
%\fi
%

\end{document}